\documentclass[12pt]{amsart}

\usepackage{amssymb}
\usepackage{amstext}
\usepackage{amsopn}
\usepackage{amsmath}
\usepackage{amsfonts}
\usepackage{amscd}

\swapnumbers

\newtheorem{theorem}[subsection]{Theorem}
\newtheorem{corollary}[subsection]{Corollary}
\newtheorem{conjecture}[subsection]{Conjecture}
\newtheorem{lemma}[subsection]{Lemma}
\newtheorem{proposition}[subsection]{Proposition}

\theoremstyle{definition}
\newtheorem{definition}[subsection]{Definition}

\newtheorem{problem}[subsection]{Problem}

\theoremstyle{remark}
\newtheorem{remark}[subsection]{Remark}

\theoremstyle{plain}
\numberwithin{equation}{subsection}

\let\cal\mathcal

\let\frak\mathfrak
\let\Bbb\mathbb

\def\BB{{\mathbf B}}

\def\FF{{\mathbf F}}

\def\LL{{\mathbf L}}

\def\QQ{{\mathbf Q}}

\def\spec{{\rm Spec}\,}

\def\cL{{\mathcal L}}
\def\cM{{\mathcal M}}
\def\cN{{\mathcal N}}
\def\cO{{\mathcal O}}
\def\cP{{\mathcal P}}

\def\cX{{\mathcal X}}

\mathchardef\alphag="7C0B
\mathchardef\betag="7C0C
\mathchardef\gammag="7C0D
\mathchardef\deltag="7C0E
\mathchardef\varepsilong="7C22
\mathchardef\varphig="7C27
\mathchardef\psig="7C20
\mathchardef\zetag="7C10
\mathchardef\epsilong="7C0F
\mathchardef\rhog="7C1A
\mathchardef\taug="7C1C
\mathchardef\upsilong="7C1D
\mathchardef\iotag="7C13
\mathchardef\thetag="7C12
\mathchardef\pig="7C19
\mathchardef\sigmag="7C1B
\mathchardef\etag="7C11
\mathchardef\omegag="7C21
\mathchardef\kappag="7C14
\mathchardef\lambdag="7C15
\mathchardef\mug="7C16
\mathchardef\xig="7C18
\mathchardef\chig="7C1F
\mathchardef\nug="7C17
\mathchardef\varthetag="7C23
\mathchardef\varpig="7C24
\mathchardef\varrhog="7C25
\mathchardef\varsigmag="7C26
\mathchardef\Omegag="7C0A
\mathchardef\Thetag="7C02
\mathchardef\Sigmag="7C06
\mathchardef\Deltag="7C01
\mathchardef\Phig="7C08
\mathchardef\Gammag="7C00
\mathchardef\Psig="7C09
\mathchardef\Lambdag="7C03
\mathchardef\Xig="7C04
\mathchardef\Pig="7C05
\mathchardef\Upsilong="7C07

\DeclareMathOperator*{\Spec}{Spec}
\def\ord{{\rm ord}}

\def\codim{{\rm codim} \,}

\def\gr{{\cal Gr}}

\begin{document}

\title[Motivic zeta functions of infinite dimensional Lie algebras]{Motivic
zeta functions of infinite dimensional Lie algebras}
\author{Marcus du Sautoy}
\address{DPMMS, Centre for Mathematical Sciences, Wilberforce Road, Cambridge CB3 0WB, U.K}
\email{dusautoy@dpmms.cam.ac.uk}
\urladdr{http://www.dpmms.cam.ac.uk/$\sim$dusautoy}
\author{Fran{\c c}ois Loeser}
\address{Centre de Math{\'e}matiques,
Ecole Polytechnique,
F-91128 Palaiseau
(UMR 7640 du CNRS), {\rm and}
Institut de Math{\'e}matiques,
Universit{\'e} P. et M. Curie, Case 82,
4 place Jussieu,
F-75252 Paris Cedex 05
(UMR 7596 du CNRS)}
\email{loeser@math.polytechnique.fr}
\urladdr{http://math.polytechnique.fr/cmat/loeser/}
\maketitle


\section{Introduction}

\subsection{}\label{1.1}

In the present paper we associate motivic
zeta functions to certain classes of infinite dimensional
Lie algebra over a field $k$ of characteristic zero.
Included in these classes are the important cases of loop
algebras, affine Kac-Moody algebras, the Virasoro algebra and Lie algebras
of Cartan type.
These zeta functions
take their values in the Grothendieck
ring of algebraic varieties over $k$ and are built by
encoding
in some manner the $k$-subalgebras of a given codimension.
This construction is done
by an adaptation of the idea of the motivic Igusa zeta function recently
introduced in \cite{Denef-Loeser-motivic Igusa}.

Let us recall
the definition
of the
Grothendieck ring ${\cal M}$ of algebraic varieties over $k$.
This is the
ring generated by symbols $[S],$ for each $S$ an algebraic variety
over $k$,
that is, a reduced separated scheme of finite type over $k$,
with the relations
\begin{enumerate}
\item[(1)] $[S]=[S^{\prime }]$ if $S$ is isomorphic to $S^{\prime }$;
\item[(2)]$[S]=[S\backslash S^{\prime }]+[S^{\prime }]$ if $S^{\prime }$ is
closed
in $S$;
\item[(3)]$[S\times S^{\prime }]=[S][S^{\prime }].$
\end{enumerate}

The idea leading to this construction comes from the analogy
with  the zeta functions capturing the subalgebra lattice of $p$%
-adic Lie algebras.
More precisely,
in \cite{GSS} the zeta function of a finite dimensional Lie algebra $L_{p}$
over ${\Bbb Z}_{p}$ was defined as
\begin{eqnarray*}
\zeta _{L_{p}}(s) &:=&\sum_{H\leq L_{p}}|L_{p}:H|^{-s} \\
&=&\sum_{n=0}^{\infty }a_{p^{n}}(L_{p})p^{-ns},
\end{eqnarray*}
where the first sum is taken over finite index ${\Bbb Z}_{p}$-subalgebras,
that is ${\Bbb Z}_{p}$-submodules which define subalgebras in $L_{p}$ and $%
a_{p^{n}}(L_{p})$ is the number of these of index $p^{n}.$ This is a
well-defined series since there are only finitely many of each for a
given finite
index. Note that since ${\Bbb Z}$ is dense in ${\Bbb Z}_{p},$ finite index $%
{\Bbb Z}_{p}$-subalgebras are the same as finite index
${\Bbb Z}$-subalgebras. The zeta function captures in some analytic manner the
lattice of subalgebras. It was proved in the same paper that, when $L_{p}$
is additively isomorphic to ${\Bbb Z}_{p}^{d}$, the zeta function is a
rational function in $p^{-s}.$ The proof depends on Denef's results on
$p$-adic definable integrals. Explicit formulae in terms of an associated
resolution of singularities are given in \cite{duSG-Abscissa} for $\zeta
_{L_{p}}(s)$ in the case that $L_{p}=L\otimes {\Bbb Z}_{p}$ for some Lie
algebra $L$ defined over ${\Bbb Z}$.

In general an infinite dimensional Lie algebra over a field $k$
of
characteristic zero contains an infinite number of subalgebras of each given
finite codimension. So it doesn't make sense to count them in a conventional
way and try to encode this counting function in a Dirichlet or Poincar{\'e}
series with integer coefficients.
Instead, in the cases we consider in the present paper, the set of
subalgebras of a fixed finite codimension form a constructible
subset of some Grassmannian, so we may asssociate to it
an element of the Grothendieck ring $\cal M$, and we can construct
the
zeta function
as a Poincar{\'e} series with coefficients in $\cal M$.

\subsection{An example~: counting
$k[[t]]$-submodules of $k[[t]]^{2}$}\label{ill}

To illustrate the thinking behind this, it is instructive to see how we
would count $k[[t]]$-submodules of finite codimension in the free $k[[t]]$%
-module $L=k[[t]]^{2}.$ Let $e_{1},e_{2}$ be a basis for $L$, then each
subalgebra $H$ of codimension $1$ can be represented by coordinates
in a $k[[t]]$-basis with respect to this standard basis. We can record these
coordinates in
a matrix whose rows are the coordinates for the basis of the subalgebra $H$.
By normalizing the matrix into upper triangular form we can choose a unique
representative
\[
\left(
\begin{array}{ll}
m_{11} & m_{12} \\
0 & m_{22}
\end{array}
\right)
\]
for each $H.$ We have two cases:

\begin{enumerate}
\item[(1)] $m_{22}\in t\left( k[[t]]\right) ^{*},$ and $m_{11}\in \left(
k[[t]]\right) ^{*}.$ To get a unique representative we can assume that $%
m_{11}=1$ and $m_{22}=t.$ By subtracting $k[[t]]$-combinations of the second
row from the first we can choose $m_{12}\in k.$ So the subalgebras of this
sort give a space of subalgebras which looks like $k.$
\item[(2)] If $m_{22}\in \left( k[[t]]\right) ^{*},$ and $m_{11}\in t\left(
k[[t]]\right) ^{*}$ then we get a unique representative of the form
\[
\left(
\begin{array}{ll}
t & 0 \\
0 & 1
\end{array}
\right) .
\]
This gives just one subalgebra of codimension $1.$
\end{enumerate}

In the case of $\left( {\Bbb Z}_{p}^{2},+\right) $ this analysis showed us
that there were $p+1$ subalgebras of index $p.$ Here we would like to say
that there are ${\Bbb A}_{k}^{1} + 1$ subalgebras, with
${\Bbb A}_{k}^{1}$ the affine line over $k$,
since, for any field $K$ containing $k$, the
$K[[t]]$-submodules of $K[[t]]^{2}$ are parametrized by $K + 1 =
{\Bbb A}_{k}^{1} (K) + 1$. For that to make sense
one is naturally led to work in the Grothendieck ring
$\cM$.
So if we
set ${\bf L}:=[{\Bbb A}_{k}^{1}]$ in $\cM$,
we can interpret our analysis of
the subalgebras in $\left( \left(
k[[t]]\right) ^{2},+\right) $ as saying there are ${\bf L}+1$ subalgebras of
codimension $1$.

On the other hand if we consider counting all vector subspaces of
codimension 1, {\it i.e.} ignoring the action of $t,$ then there are really
too many subspaces to be able to assign an element of the Grothendieck
ring
${\cal M}$. In fact, this is similar to the situation when $k={\Bbb F}_{p}$.
There are $1+p+\cdots +p^{n}$ ${\Bbb F}_{p}[[t]]$-submodules of index $%
p^{n+1}$ in $\left( {\Bbb F}_{p}[[t]]\right) ^{2}$ but an infinite number of
additive ${\Bbb F}_{p}$-subspaces of index $p^{n+1}.$

\subsection{Stability and the motivic zeta function}

The key to being able to assign an element of the Grothendieck ring to
subalgebras of a fixed codimension is a notion we call
{\em stability}.
The context in which we shall work is that of ${\Bbb N}$-filtered or ${\Bbb Z%
}$-filtered Lie algebras over a field $k$ of characteristic zero.
These are infinite dimensional Lie algebras $L$
with a filtration $L=\bigcup L_{i}$ indexed by ${\Bbb N}$ or ${\Bbb Z}$
consisting of subalgebras $L_{i}$ of finite codimension in $L_{0}.$
By a class ${\cal X}$
of subalgebras of $L$, we shall mean the data, for every
field $K$ which is a finite extension of $k$, consisting of a family ${\cal
X} (K)$
of subalgebras of $L \otimes K$.
We will sometimes, with a slight abuse, identify ${\cal X}$ with  the
union of the various sets ${\cal X} (K)$. Let $A_{n}({\cal X}) (K)$ denote
the set of
subalgebras in
${\cal X}(K)$ of codimension $n$ in  $L \otimes K$.
We say that a class of
subalgebras ${\cal X}$ of $L$ is {\em stable} if, for every $n$,
there exists
some $f(n)$ such that, for every $K$,
if $H$ belongs to $A_{n}({\cal X}) (K)$,
then $H$ contains  $L_{f(n)} \otimes K$.

In this case it is possible to show for certain classes of subalgebras
${\cal X}$ that $A_{n}({\cal X})$ is a
constructible subset of the Grassmannian ${\rm Gr}(L/L_{f(n)}).$ We can then
define the concept of a {\em motivic zeta function}:

\[
P_{L,{\cal X}}(T)
:= \sum_{n=0}^{\infty } \,[A_{n}({\cal X})] \, T^{n}
\]
where the coefficients $[A_{n}({\cal X})]$ belong to
${\cal M}.$ This zeta function
will be our interpretation of how to add the series
\[
\sum_{H\in {\cal X}}T^{\codim H},
\]
which is then an element of ${\cal M}[[T]].$

We also introduce a variant of this function, which has not been considered
previously
in the $p$-adic case, which counts subalgebras commensurable with the
subalgebra $L_{0}$ of infinite codimension in a ${\Bbb Z}$-filtered Lie
algebra $L$. A subalgebra $H$ is commensurable with $L_{0}$ if $H\cap L_{0}$
has finite codimension in $H$ and $L_{0}$. We then define the codimension of
$H$ in $L_{0}$ to be the sum of the codimension of $H\cap L_{0}$ in $H$ and $
L_{0}$.

\medskip

Our approach to proving the stability of subalgebras in a filtered Lie
algebra is a new concept we call {\em well-covered}. This means essentially
that elements in the $n$-th term of the grading can be realised in many ways
as commutators of elements further up
the grading. We prove that affine
Kac-Moody Lie algebras (which are built from loop algebras $L\otimes k[[t]]$,
with $L$ a finite dimensional Lie algebra over $k$),
the Virasoro algebra and Lie
algebras of Cartan type are well-covered and hence have an associated
motivic zeta function.

On the other hand we show that a finitely generated free Lie algebra is not
well-covered since an Engel element like $(a,b,b,...,b)$ of length $n$
cannot be realised as a commutator in more than one way. Hence we can show
that such Lie algebras don't have an associated motivic zeta function.

\subsection{Rationality of the motivic zeta function}
Once the motivic zeta function is defined, it is natural to consider
the question of its rationality.
By using results from \cite
{Denef-Loeser-Arcs} and \cite{Denef-Loeser-motivic Igusa}, we are able
to prove the
rationality of motivic zeta functions counting $k[[t]]$-subalgebras in $%
k[[t]]$-Lie algebras of the form
$L_{k}\otimes k[[t]]$ where $L_{k}$ is a finite
dimensional Lie algebra over $k$, a field of characteristic zero.

Denote by ${\cal M}_{{\rm loc}}$ the ring ${\cal M}[{\bf L}^{-1}]$ obtained
by localization and define by ${\cal M}[T]_{{\rm loc}}$ the subring of $%
{\cal M}_{{\rm loc}}[[T]]$ generated by ${\cal M}_{{\rm loc}}[T]$ and the
series $(1-{\bf L}^{a}T^{b})^{-1}$ with $a\in {\Bbb Z}$ and $b\in {\Bbb N}.$

\begin{theorem}
Let $k$ be a field of characteristic zero.

\begin{enumerate}
\item[(1)]  If $L$ is a finite dimensional free $k[[t]]$-Lie algebra of the
form $L=L_{k}\otimes _{k}k[[t]]$, with $L_{k}$ a Lie algebra over $k$,
and
${\cal X} (K)$ is the set of all $K[[t]]$-subalgebras of
$L \otimes K$,
then $P_{L,{\cal X}}(T)$ is
well-defined and is rational, belonging to ${\cal M}[T]_{{\rm loc}}$.

\item[(2)]  If $L$ is a finite dimensional $k((t))$-Lie algebra and $L_{0}$
is a choice of some $k[[t]]$-Lie subalgebra of the form $L_{0}=L_{k}\otimes
_{k}k[[t]]$, with $L_{k}$ a Lie algebra over $k$, and
${\cal X} (K)$ is the set of
all $K[[t]]$-subalgebras of $L \otimes K$
commensurable with $L_{0} \otimes K$, then $P_{L_{0},{\cal X}%
}(T)$ is well-defined and is rational, belonging to ${\cal M}[T]_{{\rm loc}}$
.
\end{enumerate}
\end{theorem}

The key to the proof of this Theorem is to express the motivic zeta function
as a motivic integral as developed by Denef and
the second author \cite{Denef-Loeser-Arcs} and \cite{Denef-Loeser-motivic
Igusa}. In fact we are naturally
led to introduce the concept of a motivic measure
on the infinite Grassmannian.

\subsection{Implications for $p$-adic zeta functions}

If $L$ is a Lie algebra over ${\Bbb Z},$ then the motivic zeta function of $%
L\otimes {\Bbb C}[[t]]$ has implications for the zeta functions of $L\otimes
{\Bbb Z}_{p}$ as we range over primes $p.$ In particular we can show, by
``taking the trace of Frobenius'' of the motivic zeta function that the
explicit formulae calculated in \cite{duSG-Abscissa} and alluded to in \ref
{1.1} for $\zeta _{L\otimes {\Bbb Z}_{p}}(s)$ are canonical and hence
independent of a choice of resolution and simplicial decomposition of a cone
which are involved in the calculation.

We also define the concept of the topological zeta function of $L$ which is
the analogue of the topological Igusa
zeta function defined in \cite
{Denef-Loeser-topological}. This can be thought of as the limit as
$p\mapsto 1$ of the explicit formulae expressing $\zeta _{L\otimes {\Bbb Z%
}_{p}}(s).$ To show that this is well-defined we show that this is the same
as taking the Euler characteristic of the associated motivic zeta function,
similarly as in
\cite{Denef-Loeser-motivic Igusa}.

\setcounter{tocdepth}{1}

\tableofcontents

{\bf Acknowledgements.} We should like to thank the London Mathematical
Society and the Institut Math{\'e}matique de Jussieu for grants which
facilitated our collaboration.  The first author would also like to thank
the Royal Society for funding in the form of a University Research
Fellowship.

\bigskip

{\bf Notations and terminology.} We shall use the notation ${\Bbb
N}=\{0,1,2,\ldots \}.$ For
$k$ a field, we denote by ${\Bbb A}_{k}^{n}$ the affine space ${\rm Spec}%
\,k[x_{1},\ldots ,x_{n}]$.
By an
algebraic variety
over $k$, we shall always mean a scheme of finite type over $k$, which
is separated and reduced.
Except explicitely stated otherwise,
all Lie algebras we consider will be assumed to be Lie algebras over $k$.
We shall consider
a measure $\widetilde{\mu }$ taking
values in ${\cal M}_{{\rm loc}}$ and a measure  $\mu $
taking
values in the completion $\widehat{{\cal M}}$ of ${\cal M}_{{\rm loc}}$.

\section{Definition of the motivic zeta function\label%
{section: definition of motivic zeta}}

We fix a field $k$. In this section we explore some infinite dimensional Lie
algebras for which one can define a motivic zeta function.

\begin{definition}
(1) An ${\Bbb N}$-filtered Lie algebra $L$ is a Lie algebra equipped with a
filtration by subspaces $L_{j}$, $j\in {\Bbb N}$, of finite codimension with
$L_{j+1}\subset L_{j}$ satisfying $(L_{j},L_{k})\subset L_{j+k}$ and $%
\bigcap_{i\in {\Bbb N}}L_{i}=0.$ We shall assume for convenience that $%
L_{0}=L.$

(2) Similarly a ${\Bbb Z}$-filtered Lie algebra comes with a filtration by
subspaces $L_{j}$, $j\in {\Bbb Z}$, with $L_{j+1}\subset L_{j}$ satisfying $
(L_{j},L_{k})\subset L_{j+k}$ where we assume that $L_{j}$ is commensurable
with $L_{0}$ and $L=\bigcup_{j\in {\Bbb Z}}L_{j}$ and $\bigcap_{i\in {\Bbb Z}%
}L_{i}=0.$
\end{definition}

The filtration defines a (formal) topology on $L$ defined by the
neighbourhood base of zero consisting of subspaces $L_{j}.$ A subspace $H$
is commensurable with $L_{0}$ if the intersection $H\cap L_{0}$ has finite
codimension in both $H$ and $L_{0}$. We define the codimension of $H$ in $%
L_{0}$ in this case as the sum of the codimension of $H\cap L_{0}$ in $H$
and $L_{0}.$ Note that in an ${\Bbb N}$-filtered Lie algebra the $L_{j}$ are
ideals but this is not necessarily the case in a ${\Bbb Z}$-filtered Lie
algebra since $(L_{-j},L_{k})\subset L_{k-j}.$

\begin{remark}
We can carry out the analysis of this section just under the assumption that
the Lie algebra $L$ comes only with a formal topology of subspaces $L_{j}$
of finite codimension without insisting that $(L_{j},L_{k})\subset L_{j+k}$.
However essentially all the cases where we can get the theory up and running
apply to the examples where the chain of subspaces is actually a filtration.
However in any particular case it may not be necessary to assume anything
this strong. It may be possible to show that subalgebras are a constructible
subset of the finite Grassmannian ${\rm Gr}(L/L_{j}).$
\end{remark}

\subsection{Examples of filtered Lie algebras}

(1) A good example of an ${\Bbb N}$-filtration on $L$ is the lower central
series. In this case $L_{0}=L_{1}$ and $L_{i}=\gamma _{i}(L)=(\gamma
_{i-1}(L),L).$ We are insisting that the filtration be of subalgebras of
finite codimension. If the first layer $L_{1}/L_{2}$ is finite dimensional
then this implies the same for the other layers $L_{i}/L_{i+1}$. If $(L,L)$
has infinite codimension then the theory won't get up and running, since
there will be too many subalgebras of codimension 1. Note that, by
definition of an ${\Bbb N}$-filtration, $L_{i}$ contains $\gamma
_{i}(L_{1}). $

(2) Another example arises if $L$ has some ${\Bbb N}$-grading, so that $%
L=\bigoplus_{j\geq 0}L(j)$ and suppose that $L(j)$ is finite dimensional
then put $L_{k}=\bigoplus_{j\geq k}L_{j}.$

(3) If $L$ has the structure of a $k[[t]]$-module then put $L_{k}=t^{k}L.$

(4) If $L$ comes with a ${\Bbb Z}$-grading $L=\bigoplus_{j\in {\Bbb Z}}L(j)$
then we can define $%
L_{k}=\bigoplus_{j\geq k}L(j)$ for $k\in {\Bbb Z}$ and consider subalgebras
commensurable with $L_{0}$.

(5) Finally if $L$ is a $k((t))$-Lie algebra with a choice of $k[[t]]$%
-subalgebra $L_{0}$ we take $L_{k}=t^{k}L_{0}$ for $k\in {\Bbb Z}$.

\subsection{Constructible sets of subalgebras in Grassmannians}

It is the subalgebras of finite codimension that are closed with respect to
the formal topology that we shall seek to count. Such subalgebras contain
then some term $L_{j}$ of the neighbourhood base and can therefore be viewed
as a point of the finite Grassmannian ${\rm Gr}(L/L_{j}).$ This follows
because if $H$ is closed, $H=\bigcap_{i\in {\Bbb N}}(H+L_{i})$ and since
$H$ has finite codimension, the chain $H+L_{i}$ must stabilize at some
point.

In the case of a closed commensurable subalgebra $H$ in a ${\Bbb Z}$%
-filtered Lie algebra, there exists some $j$ with $L_{-j}\supseteq
H\supseteq L_{j}.$ Hence these will be points of the finite Grassmannian $%
{\rm Gr}(L_{-j}/L_{j}).$ For convenience we set $X_{j}={\rm Gr}(L/L_{j})$ or
${\rm Gr}(L_{-j}/L_{j})$ according to which situation we find
ourselves.

\begin{remark}
In section \ref{section:Grassmannian} we shall introduce the concept of
infinite Grassmannians which can be thought of as the direct limit of these
finite Grassmannians. These are considered for example in Section 7.2 (i) of
\cite{Pressley-Segal}. The set of subalgebras ${\cal X}$ can then be thought
of as subsets of these infinite Grassmannians. We shall see later that in
some cases it is possible to put a motivic measure directly on an infinite
Grassmannian and express the motivic zeta functions of this section as
integrals over ${\cal X}$ with respect to this measure.
\end{remark}

\begin{remark}
Note that not all subalgebras of finite codimension will be closed with
respect to a chosen formal topology. For example in $\bigoplus_{i\in {\Bbb N}%
}{k}e_{i}$ with respect to the formal topology $L_{j}=\bigoplus_{i\geq
j}{k}e_{i},$ the subspace $$\left\langle
e_{0}-e_{1},e_{1}-e_{2},\ldots, e_{k}-e_{k+1},\ldots
\right\rangle $$ has codimension
1 but contains no term of the filtration $L_{j}.$
\end{remark}

\begin{remark}
There is a one-to-one correspondence preserving the codimension between
closed subalgebras of finite codimension of $L$ and closed subalgebras of
finite codimension of the completion of $L$ with respect to the formal
topology.
\end{remark}

\subsection{}Let $Y$ be an algebraic variety
of finite type over $k$. By the set underlying $Y$ we shall always mean
the set of
closed points of $Y$. These points correspond to rational points of
$Y$ over fields $K$ which are finite extensions of $k$.
The boolean algebra of constructible subsets
of
$Y$ is the smallest family of
subsets of $Y$ containing all  Zariski closed subsets of $Y$ and
stable by taking finite unions and complements. It follows from the
definitions
that to any
constructible subset $A$ of $Y$ one can associate a canonical
element $[A]$ of the Grothendieck
ring ${\cal M}$.

We shall now view the finite Grassmannian $X_{j}={\rm Gr}(L/L_{j})$ or
${\rm Gr}(L_{-j}/L_{j})$ as an
algebraic variety
of finite type over $k$. Hence points of $X_{j}$ will be linear spaces
defined over some field $K$ which is a finite extension of $k$.

By a class ${\cal X}$
of subalgebras of $L$, we shall mean the data, for every
field $K$ which is a finite extension of $k$, consisting of a family ${\cal
X} (K)$
of subalgebras of $L \otimes K$.
For every $l,n\in {\Bbb N}$, and
every
field $K$ which is a finite extension of $k$, we define the subsets$:$%
\begin{eqnarray*}
A_{n}({\cal X})(K) &:=&\left\{ H\in {\cal X}(K):
{\rm codim}_{L \otimes K}
\,H=n\right\} \\
A_{l,n}({\cal X}) (K)
&:=&\left\{ H\in A_{n}({\cal X})(K) :H\supseteq L_{l}\otimes K\right\}.
\end{eqnarray*}
Hence $A_{l,n}({\cal X})(K)$ can be thought of then as a subset of
$X_{l}(K)$, and the union of these subsets when varying $K$ is a
subset
$A_{l,n}({\cal X})$ of the finite Grassmannian
$X_{l}$.

\begin{definition}
We call ${\cal X}$ a {\em constructible class of subalgebras} if $A_{l,n}(%
{\cal X})$ is a constructible subset of $X_{l}$, for every
$l,n\in {\Bbb N}$.
\end{definition}

\subsection{Examples of constructible classes of subalgebras}

For $L$  an ${\Bbb N}$-filtered Lie
algebra, we define ${\cal X}^{\leq }$ as the class of closed
subalgebras, {\it i.e.} ${\cal X}^{\leq } (K)$ is the set of all
closed subalgebras of $L \otimes K$, and
${\cal X}^{\triangleleft }$
as the class of closed
subalgebras which are ideals,
{\it i.e.} ${\cal X}^{\triangleleft} (K)$ is the set of all
closed subalgebras of $L \otimes K$ which are ideals of
$L \otimes K$.

If $L$ is a ${\Bbb Z}$-filtered Lie algebra, the corresponding
classes of commensurable subalgebras are defined as follows: for $*\in
\left\{ \leq ,\triangleleft
\right\}$,
${\cal X}_{0}^{*} (K)$ is the set of
subalgebras $H$ in  ${\cal X}^{*} (K)$ which are
commensurable with
$L_{0} \otimes K$.

\begin{lemma}
\label{n-filtered implies constructible}Let $L$ be an ${\Bbb N}$-filtered
Lie algebra then for $*\in \left\{ \leq ,\triangleleft \right\} ,$ ${\cal X}%
^{*}$ is a constructible class of subalgebras.
\end{lemma}

\begin{proof}Choose a basis $e_{1},\ldots ,e_{r}$ for $L/L_{l}.$ Let $\beta
: {\Bbb A}_k ^{2r}\rightarrow {\Bbb A}_k ^{r}$ denote the bilinear form
defining the Lie bracket
in $L/L_{l}$ with respect to this basis. Let ${\rm Tr}_{r, k}$ denote the
$k$-variety of
upper triangular matrices of order $r$. Hence the $K$-points of
${\rm Tr}_{r, k}$ are just the upper triangular matrices of order $r$
with coefficients in $K$.
We have a morphism
\[
h:{\rm Tr}_{r, k} \longrightarrow {\rm Gr}(L/L_{l})
\]
which takes the matrix $\left( m_{ij}\right) $ to the subspace spanned by
\[
\left\{ m_{11}e_{1}+\cdots +m_{1r}e_{r},\ldots ,m_{rr}e_{r}\right\} .
\]
Let ${\cal M}_{l,n}^{*}$ denote the inverse image of $A_{l,n}({\cal X}^{*})$
under this map. Since there is a finite partition of ${\rm Gr}(L/L_{l})$
into locally closed subvarieties over which $h$ induces morphisms which are
locally trivial fibrations for the Zariski topology, it suffices to show
that ${\cal M}_{l,n}^{*}$ is a constructible subset of ${\rm Tr}_{r,k}$.

\begin{enumerate}
\item[(1)] a matrix $\left( m_{ij}\right) $ in
${\rm Tr}_{r, k} (K)$ defines a subspace of codimension $n$ in $L
\otimes K$
if and only if the diagonal $\left( m_{11}, \ldots, m_{rr}\right) $ contains
exactly $n$ zero entries.

\item[(2$^{\leq }$)] a matrix $\left( m_{ij}\right) $ in
${\rm Tr}_{r, k} (K)$
defines a subalgebra of $L \otimes K$
if and
only if for each $1\leq i<j\leq r$ there exist $Y_{ij}^{1}, \ldots,
Y_{ij}^{r}\in K$ such that
\[
\beta ({\bf m}_{i},{\bf m}_{j})=\sum_{k=1}^{r}Y_{ij}^{k}{\bf m}_{k}
\]
where ${\bf m}_{i}$ denotes the $i$th row of $\left( m_{ij}\right) ;$

\item[(2$^{\triangleleft }$)] a matrix $\left( m_{ij}\right) $ in
${\rm Tr}_{r, k} (K)$
defines an ideal of $L \otimes K$ if
and only if for each $1\leq i,j\leq r$ there exist $Y_{ij}^{1}, \ldots,
Y_{ij}^{r}\in K$ such that
\[
\beta ({\bf e}_{i},{\bf m}_{j})=\sum_{k=1}^{r}Y_{ij}^{k}{\bf m}_{k}
\]
where ${\bf e}_{i}$ denotes the $r$-tuple with 1 in the $i$-th entry and zeros
elsewhere.
\end{enumerate}

The two conditions on a matrix $\left( m_{ij}\right) $ in (1) and (2$^{*}$)
define ${\cal M}_{l,n}^{*}$ as a constructible subset of
${\rm Tr}_{r, k}$.
\end{proof}

If $L$ is a finite dimensional $k[[t]]$-Lie algebra then we can also
define the
classes of subalgebras
${\cal X}_{t}^{\leq }$
with
${\cal X}_{t}^{\leq } (K)$ the set of closed subalgebras of $L \otimes K$
which are
$K [[t]]$-submodules of $L \otimes K$
and
${\cal X}_{t}^{\triangleleft }$
defined by
${\cal X}_{t}^{\triangleleft } (K)
=
{\cal X}_{t}^{\leq }(K)\cap {\cal X}^{\triangleleft }(K)$.

\begin{lemma}
Let $L$ be a finite dimensional $k[[t]]$-Lie algebra whose filtration
consists of the ideals $L_{j}=t^{j}L$. Then for $*\in \left\{ \leq
,\triangleleft \right\} ,$ ${\cal X}_{t}^{*}$ are constructible classes of
subalgebras.
\end{lemma}

\begin{proof}We just have to add the constructible condition:
\begin{enumerate}
\item[(3)] a matrix $\left( m_{ij}\right) $ in
${\rm Tr}_{r, k} (K)$
defines a $K[[t]]$-submodule of $L \otimes K$
if and
only if for each $1\leq i\leq r$ there exist $Y_{i}^{1}, \ldots,
Y_{i}^{r}\in K$ such that
\[
t{\bf m}_{i}=\sum_{k=1}^{r}Y_{i}^{k}{\bf m}_{k}.
\]
\end{enumerate}
\end{proof}

We shall see later a better reason for these subalgebras being constructible
when we realise ${\cal X}_{t}^{*}$ as a semi-algebraic subset in a suitable
arc space.

\begin{lemma}
Suppose that $L$ is a ${\Bbb Z}$-filtered Lie algebra. Then ${\cal X}_{0}^{%
\leq }$ is a constructible set of subalgebras.
\end{lemma}

\begin{proof}We want to put a condition on a subspace of $L_{-j}/L_{j}$ that
it is a subalgebra of codimension $n$ in $L_{0}.$ The trouble is now that $%
L_{j}$ is not necessarily an ideal in $L_{-j}.$ However we use the fact that
$(L_{-j},L_{2j})\subset L_{j}.$ Let $e_{0},e_{-1},\ldots ,e_{-s_{1}}$ be a
basis for $L_{-2j}/L_{0}$ with $e_{0},e_{-1},\ldots ,e_{-s}$ a basis for $
L_{-j}/L_{0}$, and $e_{1},\ldots ,e_{r_{1}}$ be a basis for $L_{0}/L_{2j}$
with $e_{1},\ldots ,e_{r}$ a basis for $L_{0}/L_{j}.$ Let $\beta
: {\Bbb A}_k ^{s+r_{1}+1}\times {\Bbb A}_k ^{s+r_{1}+1}\rightarrow {\Bbb
A}_k ^{s_{1}+r+1}$ be the bilinear
map defining the Lie bracket from $L_{-j}/L_{2j}\rightarrow L_{-2j}/L_{j}$
with respect to these choices of basis.

There is still a morphism
\[
{\rm Tr}_{s+r+1, k}
\rightarrow {\rm Gr}(L_{-j}/L_{j})
\]
which takes a matrix $\left( m_{ij}\right) $, $-s\leq i\leq j\leq r$, to the
subspace spanned by
\[
\left\{ m_{-s,-s}e_{-s}+\cdots +m_{-s,r}e_{r}, \ldots, m_{rr}e_{r}\right\} .
\]

(1) a matrix $\left( m_{ij}\right) $ in
${\rm Tr}_{s+r+1, k} (K)$ defines a subspace of codimension $n$
in $L_{0} \otimes K$ if and only if the diago\-nal $\left( m_{11}, \ldots,
m_{rr}\right) $ contains exactly $n_{1}$ zero entries and the diagonal $%
(m_{-s,-s}, \ldots, m_{00})$ contains $n_{2}$ non zero entries and $%
n_{1}+n_{2}=n$.

(2$^{\leq }$) the subspace $H=\left\langle m_{-s,-s}e_{-s}+\cdots
+m_{-s,r}e_{r},\ldots ,m_{rr}e_{r}\right\rangle +L_{j} \otimes K$
defines a subalgebra of $L \otimes K$
if and only if, for $-s\leq i<j\leq r$ and
for all $\lambda =\left( \lambda
_{r+1},\ldots ,\lambda _{r_{1}}\right) ,\mu =\left( \mu _{r+1},\ldots ,\mu
_{r_{1}}\right) \in K^{r_{1}-r}$ there exist $Y_{ij}^{-s},\ldots
,Y_{ij}^{r}\in K$ such that
\[
\beta ({\bf m}_{i},\lambda ,{\bf m}_{j},\mu )=\left( 0,\ldots
,0,\sum_{k=-s}^{r}Y_{ij}^{k}{\bf m}_{k}\right)
\]
where there are $s_{1}-s$ zeros. The point is that we are going to be
guaranteed $(H,L_{2j})\subset H$ so we just need to check that $H$ is a
subalgebra modulo $L_{2j}$ which is finite dimensional and therefore leads
to a constructible subset. \end{proof}

Proving that ideals in a ${\Bbb Z}$-filtered Lie algebra define a
constructible set is slightly more problematic since we need to check that
the action of all the $L_{-j}$ stabilise the candidate ideal $H.$ However,
if $L$ is finitely generated we need only check these finite number of
generators:

\begin{lemma}
Suppose that $L$ is a finitely generated ${\Bbb Z}$-filtered Lie algebra.
Then ${\cal X}_{0}^{\triangleleft }$ is a constructible set of subalgebras.
\end{lemma}

\begin{proof}Let $f_{1}, \ldots, f_{d}$ be a finite set of
generators. Then, for any fixed $j$, there exists some $N(j)$ and $M(j)$ such
that $(f_{i},L_{N(j)})\subset L_{j}$ and $(f_{i},L_{-j})\subset L_{-M(j)}.$
So we can carry out the same analysis essentially as the previous lemma. Let
$e_{0},e_{-1},\ldots,e_{-s_{1}}$ be a basis for $L_{-M(j)}/L_{0}$ with $
e_{0},e_{-1},\ldots,e_{-s}$ a basis for $L_{-j}/L_{0}$, and $%
e_{1},\ldots,e_{r_{1}} $ be a basis for $L_{0}/L_{N(j)}$ with $
e_{1},\ldots,e_{r}$ a basis for $L_{0}/L_{j}.$ Let
$\varphi_{i}: {\Bbb A}_k^{s+r_{1}+1}\rightarrow {\Bbb A}_k^{s_{1}+r+1}$
define the action of $%
f_{i}$ taking $L_{-j}/L_{N(j)}\rightarrow L_{-M(j)}/L_{j}$ with respect to
these choices of basis. Then

(2$^{\triangleleft }$) the subspace $H=\left\langle m_{-s,-s}e_{-s}+\cdots
+m_{-s,r}e_{r},\ldots,m_{rr}e_{r}\right\rangle +L_{j} \otimes K$
defines
an ideal of $L \otimes K$ if
and only if it is a subalgebra of $L \otimes K$
and for $-s\leq j\leq r$ and $i=1,\ldots,d$
for all $\lambda =\left( \lambda_{r+1},\ldots,\lambda_{r_{1}}\right) \in
K^{r_{1}-r}$ there exist $Y_{ij}^{-s},\ldots,Y_{ij}^{r}\in K$ such that
\[
\varphi_{i}({\bf m}_{i},\lambda )=\left( 0,\ldots,0,\sum_{k=-s}^{r}Y_{ij}^{k}%
{\bf m}_{k}\right)
\]
where there are $s_{1}-s$ zeros. \end{proof}

Let $L$ be a finite dimensional $k((t))$-Lie algebra and $L_{0}$ be some
choice of a $k[[t]]$-submodule of the same dimension. Set $L_{j}=t^{-j}L_{0}$
for $j\in {\Bbb Z}.$ Note that in this case there are no commensurable
ideals unless $L$ is abelian. The reason for this is as follows. Let
$e_{1},\ldots ,e_{d}$ be a basis for the $
k[[t]]$-submodule $L_{0}.$ Then commensurable ideals must be contained in
the centre of $L$ since if $H\leq L_{-j}$ is an ideal and $\left(
e_{i},h\right) =x\in L_{k}\backslash L_{k+1}$ then $\left(
t^{-k-j-1}e_{i},h\right) =t^{-k-j-1}x\notin H.$ But the centre $Z(L)$ is a $%
k((t))$-submodule. So if there exists a commensurable ideal $H$ then there
exists some $L_{j}\leq H\leq Z(L).$ Hence $Z(L)=L.$ So the case of
commensurable ideals in finite dimensional $k((t))$-Lie algebras is not an
interesting one.

In the same manner as above we can prove:

\begin{lemma}
Suppose that $L$ is a finite dimensional $k((t))$-Lie algebra. Then
${\cal X}_{t,0}^{\leq }$, defined by
\[
{\cal X}_{t,0}^{\leq } (K)=\left\{ H\in {\cal X}_{0}^{\leq } (K):
H\text{ is a }K[[t]]%
\text{-submodule of } L \otimes K\right\},
\]
is a constructible class of subalgebras.
\end{lemma}

\begin{proof}Again we just have to add the constructible condition:

(3) a matrix $\left( m_{ij}\right) $ in ${\rm Tr}_{s+r+1, k} (K)$
defines a $K[[t]]$-submodule if and
only if, for each $-s\leq i\leq r$, there exist $Y_{i}^{-s},\ldots,Y_{i}^{r}%
\in K$ such that
\[
t{\bf m}_{i}=\sum_{k=-s}^{r}Y_{i}^{k}{\bf m}_{k}.
\]
\end{proof}

\subsection{Stable classes of subalgebras and motivic zeta functions}

If ${\cal X}$ is a constructible class of subalgebras then $A_{l,n}({\cal X}%
) $ defines an element $\left[ A_{l,n}({\cal X})\right] $ of the
Grothendieck ring ${\cal M}$. The only way that this can have a limit as $
l\rightarrow \infty $ is that the series stabilize at some point.

\begin{definition}
We shall call a constructible class of subalgebras ${\cal X}$ a {\em stable
class of subalgebras }if for every $n$ there
exists some $l$ such that $A_{l,n}=A_{k,n}$ for all $k\geq l$.
When this is the case we set
$\left[ A_{n}\right] =\left[ A_{l,n}\right] $ and we
define the {\em motivic zeta function} of
$L$ and ${\cal X}$ as
\[
P_{L,{\cal X}}(T)=\sum_{n=0}^{\infty }\left[ A_{n}\right] T^{n}.
\]
\end{definition}

Notice we can now say why it isn't sensible to count ${k}$-subalgebras of
finite codimension in $k[[t]]^{2}$ with trivial Lie structure since they are
not a stable system of subalgebras. As we said above, there are too many
such subalgebras. In fact ${\cal X}^{\leq }$ (respectively ${\cal X}
^{\triangleleft }$) will not be stable for all ${\Bbb N}$-filtered or ${\Bbb %
Z}$-filtered Lie algebras $L$ for which $\overline{(H,H)}$ (respectively $%
\overline{(L_{0},H)})$ has infinite codimension for any closed subalgebra $H$
(respectively ideal $H)$ of finite codimension.

\subsection{Stability of ideals}

If we only want to count ideals then we can prove stability for all ${\Bbb N}
$-filtered Lie algebras, for which the closure of terms in the lower central
series of $L_{1}$ have finite codimension:

\begin{theorem}
Suppose $L=\bigcup_{i\in {\Bbb N}}L_{i}$ is an ${\Bbb N}$-filtered Lie
algebra and $\overline{\gamma _{i}(L_{1})}$ has finite codimension for all $i
$ (where $\overline{\gamma _{i}(L)}$ is the closure of the $i$th term of the
lower central series with respect to the formal topology defined by $L_{i})$%
. Then ${\cal X}^{\triangleleft }(L)$ is a stable class of subalgebras.
\end{theorem}

\begin{proof}Suppose $H$ is a closed ideal of codimension $n$ in $L.$
Let us show that each ideal $H\cap
\bigcup_{i\geq 1}L_{i}$ of codimension $n$ in $\bigcup_{i\geq 1}L_{i}$
contains some fixed term of the filtration. We can
assume $L_{0}=L_{1}$.
Let us  consider $H_{i}=\left( H\cap
\gamma _{i}(L)\right) \gamma _{i+1}(L)$ where $\gamma _{i}(L)$ is the lower
central series of $L$ (which might be distinct from the closure of the lower
central series with respect to the filtration $L_{i})$. The ideal
$H$ is closed and
hence contains some term of the filtration $L_{m}\supset \gamma _{m}(L).$
There must be some $1\leq i\leq n+1$ such that $H_{i}=$ $\gamma _{i}(L)$
else the ideal would have codimension greater than $n.$ Hence $H+\gamma
_{i+1}(L)\supset \gamma _{i}(L).$ So $H\cap \gamma _{i}(L)\equiv \gamma
_{i}(L)\mod \gamma _{i+1}(L).$ Now $H\cap \gamma _{i+1}(L)\supset (L,H\cap
\gamma _{i}(L))$ since $H$ is an ideal. Hence $H\cap \gamma _{i+1}(L)\equiv
\gamma _{i+1}(L)\mod \gamma _{i+2}(L)$ and therefore $H+\gamma
_{i+2}(L)\supset H+\gamma _{i+1}(L)\supset \gamma _{i}(L).$ Continuing this
analysis shows that $H\cap \gamma _{i+e}(L)\equiv \gamma _{i+e}(L)\mod %
\gamma _{i+e+1}(L)$ for all $e$ hence $H+\gamma _{i+e+1}(L)\supset H+\gamma
_{i+e}(L)\supset \gamma _{i}(L).$ But $H+\gamma _{i+e+1}(L)=H$ for some $e$
since $H$ is a closed subalgebra of finite codimension. This implies that $%
H\supset \gamma _{n+1}(L)$. Since $H$ is closed, $H\supset \overline{\gamma
_{n+1}(L)}\supset L_{f(n)}$.
Similarly, if $K$ is a finite extension of $k$, every closed ideal of
codimension $n$ in $L \otimes K$ contains $L_{f(n)} \otimes K$,
whence ${\cal X}^{\triangleleft }(L)$
is a stable class of subalgebras. \end{proof}

\begin{corollary}
Suppose $L=\bigcup_{i\in {\Bbb N}}L_{i}$ is an ${\Bbb N}$-filtered Lie
algebra and $L_{0}=L_{1}.$ Then ${\cal X}^{\triangleleft }(L)$ is a stable
class of subalgebras if and only if $\overline{\gamma _{i}(L)}$ has finite
codimension for all $i.$
\end{corollary}

Note that $\overline{\gamma _{i}(L)}$ has finite codimension for all $i$ if
and only if $\overline{(L_{1},L_{1})}$ has finite codimension.

\begin{remark}
There is a case missing in the above analysis: when $L_{0}\neq L_{1},$ $%
\overline{(L_{0},H)}$ has finite codimension for all closed ideals $H$ of
finite codimension but $\overline{(L_{1},L_{1})}$ has infinite codimension.
\end{remark}

\subsection{Stability of subalgebras}

We now describe several important categories of Lie algebras with stable
classes of subalgebras. Note that for two classes of constructible
subalgebras ${\cal X}^{\natural }(L)\subset {\cal X}^{\flat }(L)$ once we
have proved that ${\cal X}^{\flat }(L)$ is stable, then ${\cal X}^{\natural
}(L)$ will also be stable.

\begin{theorem}
\label{k[[t]] algebra are stable}(1) Let $L$ be a finitely dimensional Lie
algebra over $k[[t]]$ and define the formal topology on $L$ via $
L_{l}=t^{l}L.$ Then ${\cal X}_{t}^{*}$ is a stable constructible class of
subalgebras for $*\in \left\{ \leq ,\triangleleft \right\} $.

(2) Suppose that $L$ is a finite dimensional $k((t))$-Lie algebra. Then $%
{\cal X}_{t,0}^{*}$ for $*\in \left\{ \leq ,\triangleleft \right\} $ is a
stable class of subalgebras.
\end{theorem}

\begin{proof}This follows from the fact that a $K[[t]]$-submodule $H$ of
codimension $n$ in $L \otimes K$
must contain $L_{n} \otimes K$. In (2) we apply this to
$H\cap (L_{0}\otimes K)$
which must have codimension bounded by $n$. Recall that the case ${\cal X}
_{t,0}^{\triangleleft }$ either is empty or else $L$ is abelian and ${\cal X}%
_{t,0}^{\triangleleft }={\cal X}_{t,0}^{\leq }.$ \end{proof}

\subsection{Loop algebras}\label{loop}

Examples of Lie algebras that are finite-dimensional as ${\Bbb C}[[t]]$-Lie
algebras arise as the positive parts of loop algebras. These are defined as
the complexified algebra of smooth maps from the circle $S^{1}$ (or more
generally some compact complex manifold) into a finite dimensional ${\Bbb C}$
-Lie algebra ${\frak g}_{{\Bbb C}},$ and denoted $L{\frak g}_{{\Bbb C}}$. We
can decompose $L{\frak g}_{{\Bbb C}}$ into its Fourier components:
\[
L{\frak g}_{{\Bbb C}}=\bigoplus_{k\in {\Bbb Z}}{\frak g}_{{\Bbb C}}\cdot
t^{k}.
\]
(One uses the notation $L{\frak g}$ to denote the set of loops, {\it \ i.e.}
the set of smooth maps before we have complexified.) This is the
decomposition into eigenspaces of the action of the circle which rotates the
loops bodily. We define the positive part of the loop algebra to be
\[
L_{0}{\frak g}_{{\Bbb C}}=\bigoplus_{k\geq 0}{\frak g}_{{\Bbb C}}\cdot
t^{k}.
\]
This is then isomorphic to the ${\Bbb C}[[t]]$-Lie algebra ${\frak g}_{{\Bbb %
C}}\otimes {\Bbb C}[[t]].$

These examples are then graded Lie algebras, unlike a general ${\Bbb C}[[t]]
$-Lie algebra arising from a bilinear form defined over ${\Bbb C}[[t]]$ on a
finite dimensional ${\Bbb C}[[t]]$-module, which are not in general graded
but just ${\Bbb N}$-filtered.

Note that $L_{0}{\frak g}_{{\Bbb C}}\otimes {\Bbb C}((t))$ is the ${\Bbb C}$%
-subalgebra generated by the rational loops sometimes denoted $L_{{\rm rat}}%
{\frak g}_{{\Bbb C}}$ (see section 3.5 of \cite{Pressley-Segal}). Also it
should be pointed out that $L_{{\rm rat}}{\frak g}_{{\Bbb C}}$ is
${\Bbb Z}$-filtered
with $L_{{\rm rat}}{\frak g}_{{\Bbb C}}=\bigcup_{i\in {\Bbb Z}%
}t^{i}L_{0}{\frak g}_{{\Bbb C}}$ whilst the Hilbert space of all loops $L%
{\frak g}_{{\Bbb C}}$ is something bigger containing $L_{{\rm rat}}{\frak g}%
_{{\Bbb C}}$ as a dense subspace.

There is another important ${\Bbb C}$-subalgebra $L_{{\rm pol}}{\frak g}_{%
{\Bbb C}}={\frak g}_{{\Bbb C}}\otimes {\Bbb C}[z,z^{-1}]$ generated by all
the loops which are given by Laurent polynomials in $z$ and $z^{-1}.$ This
is the same as the Lie algebra of polynomial maps from ${\Bbb C}^{*}$ into $%
{\frak g}$ - ${\Bbb C}^{*}$ is the complexification of the circle (see 2.2
of \cite{Pressley-Segal}). In the case that ${\frak g}_{{\Bbb C}}$ is a
simple finite dimensional ${\Bbb C}$-Lie algebra, the subalgebras of
polynomial loops are what the Kac-Moody algebras are built from. Kac-Moody
Lie algebras arise out of the existence of interesting central extensions
that these loop algebras have. We shall see that affine Kac-Moody algebras
are examples of infinite dimensional Lie algebras for whom the class of all
subalgebras is a
stable class.

\subsection{Kac-Moody Lie algebras}

A Kac-Moody Lie algebra over a field $k$ of characteristic zero is defined
in the first instance via a presentation
associated to a generalized Cartan matrix (see \cite{Kac}). The affine
Kac-Moody algebras are those whose Cartan matrices have corank 1. They can
then be realised as algebras whose derived group is the universal central
extension of the Lie algebra $L_{{\rm pol}}{\frak g}_k={\frak g}_k\otimes
k[t,t^{-1}]$ of polynomial loops from $k^{*}$ into a
simple finite dimensional Lie algebra ${\frak g}_{k}$. More generally
if ${\frak g}_{k}={\frak g}_{1,k}\oplus \cdots \oplus {\frak g}%
_{q,k}$ is a decomposition of ${\frak g}_{k}$ into simple
factors then the product of algebras $k\cdot \tilde{\oplus}\tilde{L}{\frak g%
}_{j,k}$ is an affine Kac-Moody Lie algebra where $k\cdot d\tilde{
\oplus}\tilde{L}{\frak g}_{j,k}$ means the semi-direct product where
the factor $k\cdot d$ is generated by the derivation $d=t\frac{d}{dt}$ of
$\tilde{L}{\frak g}_{j,k}$ and $\tilde{L}{\frak g}%
_{j,k}$ is the universal central extension of the
polynomial loop algebra $L_{{\rm pol}}{\frak g}_{j,k}$ which looks
additively like $k\cdot c\oplus L_{{\rm pol}}{\frak g}_{j,k}.$ Not
all are covered by this construction. The other affine Kac-Moody Lie
algebras arise by taking the same construction associated to the {\em %
twisted loop algebras}.

These are defined for a choice of outer automorphism $\alpha $ of ${\frak g}%
_{k}$ of finite order $m$ and for fields $k$ containing primitive $m$-th
roots of unity. One replaces $L_{{\rm pol}}{\frak g}_{%
k}$ by its subalgebra $L_{(\alpha )}{\frak g}_{k}$ consisting
of the loops $f\in L_{{\rm pol}}{\frak g}_{k}$ which are equivariant:
\[
f(\varepsilon ^{-1}t)=\alpha (f(t))
\]
where $\varepsilon $ is a primitive $m$-th root of unity. The automorphism
gives rise to a decomposition of the finite dimensional Lie algebra:
\[
{\frak g}_{k}=\bigoplus_{j\in {\Bbb Z}/m{\Bbb Z}}{\frak g}_{j},
\]
where ${\frak g}_{j}$ is the eigenspace of $\alpha $ for the eigenvalue $%
\varepsilon ^{j}.$ Conversely any ${\Bbb Z}/m{\Bbb Z}$-grading of ${\frak g}%
_{k}$ arises like this since the linear transformation of ${\frak g}_{k}$
given by multiplying the vectors of ${\frak g}_{j}$ by $
\varepsilon ^{j}$ is an automorphism $\alpha $ of ${\frak g}_{k}$
which has order $m.$ The twisted loop algebra then has the following
description as a ${\Bbb Z}$-graded Lie algebra:
\[
L_{(\alpha )}{\frak g}_{k}=\bigoplus_{n\in {\Bbb Z}}{\frak
  g}_{n \, {\rm mod} %
m}\otimes t^{n}.
\]
(This is sometimes referred to as a loop algebra by algebraists because of
the cyclic filtration and $L_{{\rm pol}}{\frak g}_{k}$ as the loop
algebra with trivial $C_{m}$-grading. This seems to be a little confusing
and perhaps arose out of a misunderstanding that the original name loop
arose from the setting over ${\Bbb C}$ of maps from the circle (hence the
name loop) rather
than a graded Lie algebra whose filtration somehow comes from `looping' the
finite graded Lie algebra ${\frak g}.)$

For more details see section 5.3 of \cite{Pressley-Segal}. Note that in our
context there is little difference between taking polynomial loops or the
rational loop space because polynomial loops are dense in the rational loops
and the lattices of closed subalgebras of finite codimension will be in one
to one correspondence.

\begin{theorem}
If $L$ is an affine Kac-Moody Lie algebra with no infinite dimensional
abelian quotients and $L_{0}$ is its positive part, then for $*\in \left\{
\leq ,\triangleleft \right\} ,$ ${\cal X}^{*}(L_{0})$ and ${\cal X}
_{0}^{*}(L)$ are stable classes of subalgebras.
\end{theorem}

\begin{proof}Let $H$ be a closed subalgebra of $L_{0}$ of codimension $n.$
Suppose $%
{\frak g}_{k}={\frak g}_{1,k}\oplus \cdots \oplus {\frak g}_{q,%
k}$ is the decomposition of the underlying simple Lie algebra. Then
\[
L_{0}=\bigoplus_{i=1}^{q}L_{(\alpha_{i}),0}{\frak g}_{i,k}\oplus
k\cdot c_{i}\oplus k\cdot d_{i}
\]
where $c_{i}$ is central, $d_{i}$ is a derivation and $\alpha_{i}$ is an
automorphism of ${\frak g}_{k}$ which has order $m$ ($k$ is assumed to have
primitive $m$-th roots of unity) and
\[
L_{(\alpha ),r}{\frak g}_{k}=\bigoplus_{n\geq r}{\frak g}_{n \,
  {\rm mod}
m}\otimes t^{n}.
\]

We can restrict our attention to showing that a subalgebra $H_{i}$ of
codimension $n_{i}$ in $L_{(\alpha_{i}),0}{\frak g}_{i,k}$ must
contain some fixed term $L_{(\alpha_{i}),f(n_{i})}{\frak g}_{i,k}$ of
the filtration on $L_{(\alpha_{i}),0}{\frak g}_{i,k}$ depending only
on $n_{i}.$ For convenience we drop the subscript $i$ and since we are only
considering one twisted algebra at a time we put $L_{r}=L_{(\alpha ),r}%
{\frak g}_{k}.$

The assumption that $L$ has no infinite abelian sections means that ${\frak g%
}_{k}$ is not abelian and hence ${\frak g}_{k}=\left( {\frak g}
_{k},{\frak g}_{k}\right) $ is perfect. This implies in turn
that $\left( L_{im},L_{jm}\right) =L_{(i+j)m}.$ Let
\[
H_{i}=\left( H\cap L_{im}\right) +L_{(i+1)m}/L_{(i+1)m}.
\]
The codimension of $H$ is then equal to the sum of the codimensions of $%
H_{i} $ in $L_{im}/L_{(i+1)m}.$ Therefore $H_{i}\neq $ $L_{im}/L_{(i+1)m}$
for at most $n$ values of $i.$ Therefore consider $N\geq 2n,$ there must be
some value of $i$ with $1\leq i\leq N/2$ such that both $%
H_{i}=L_{im}/L_{(i+1)m}$ and $H_{N-i}=L_{(N-i)m}/L_{(N-i+1)m}.$ Hence $%
H_{N}\supset \left( H_{i},H_{N-i}\right) =L_{Nm}/L_{(N+1)m}$ for all $N.$
Thus $H\supset L_{2nm}. $\end{proof}

Note that an affine Kac-Moody Lie algebra is finitely generated hence ${\cal %
X}_{0}^{\triangleleft }(L)$ is constructible.

The key to this proof is that in affine Kac-Moody Lie algebras, each
element can be
realised as a commutator in many different ways.

In fact we can see that the argument applies to the following filtered Lie
algebras:

\begin{definition}
Call an ${\Bbb N}$-filtered Lie algebra {\em well-covered} if it has the
property that for each $n$ there exists some $f(n)$ such that $%
L_{f(n)}=(L_{i},L_{f(n)-i})$ for $n+1$ values of $i<f(n)/2.$ Call a ${\Bbb Z}
$-filtered Lie algebra well-covered if its positive part is well-covered.
\end{definition}

The proof of the theorem above can then easily be adapted to prove:

\begin{theorem}
(1) If the  ${\Bbb N}$-filtered Lie algebra $L$
is well-covered then ${\cal X}%
^{*}(L)$ for $*\in \left\{ \leq ,\triangleleft \right\} $ are stable
classes of subalgebras.

(2) If the ${\Bbb Z}$-filtered Lie algebra $L$
is well-covered then ${\cal X}_{0}^{%
\leq }(L)$ is a stable class of subalgebras.
\end{theorem}

The affine Kac-Moody Lie algebras (or rather the underlying loop algebras)
make up the bulk of the classification of the simple ${\Bbb Z}$-graded Lie
algebras over an algebraically closed field of characteristic zero whose
graded pieces have finite growth. According to the classification of such
Lie algebras proved by Mathieu \cite{Mathieu} the remaining Lie algebras are
those of Cartan type and the Virasoro algebra.

\subsection{Virasoro algebra}\label{vir}

The Virasoro algebra is defined for any field $k$ of characteristic zero as
the unique central extension of the Witt algebra. The Witt algebra ${\frak
d}$ is defined by:
\[
{\frak d}:={\rm Der}k[t,t^{-1}]=\bigoplus_{j\in {\Bbb Z}}k\cdot d_{j}
\]
where $d_j=-t^{j+1}\frac{d}{dt}$ with the following commutation relations:
\[
\left[ d_{i},d_{j}\right] =(i-j)d_{i+j}.
\]

The Lie algebra ${\frak d}$ has a unique (up to isomorphism) non-trivial
central extension by a 1-dimensional centre, $k\cdot c$ say, called the
Virasoro algebra ${\rm Vir}$, which is defined by the following commutation
relations:
\[
\left[ d_{i},d_{j}\right] =(i-j)d_{i+j}+(i^{3}-i)\delta _{i,-j}c\text{ where
}i,j\in {\Bbb Z}.
\]

The Virasoro algebra, denoted ${\rm Vir}$, plays an important role in the
representation theory of affine Kac-Moody Lie algebra and in quantum field
theory. It
is a ${\Bbb Z}$-graded subalgebra in ${\rm Der}\tilde{L}{\frak g}_{k%
}, $ the group of derivations of the algebra $\tilde{L}{\frak g}_{k},$.
When $k={\Bbb C}$,
it is the complexified Lie algebra of smooth vector fields on the circle.
The positive part of this Lie algebra ${\frak d}^{+}:=\bigoplus_{j\geq 1}%
k\cdot d_{j}$ also arises as the graded Lie algebra of the so-called
Nottingham group over $k.$ This is defined as the group whose
underlying set is the set of formal power series $tk[[t]]$ and whose
group operation is substitution of power series. It is a certain subgroup of
the automorphism group of $k[[t]],$ called the wild automorphisms.

\subsection{Cartan subalgebras}

The Cartan algebras are defined as follows. Let $n\geq 1$ and let $\mathbf{W}%
_{n}$ be the algebra of derivations of the polynomial ring $R=k[X_{1},\ldots
,X_{n}].$ So $\mathbf{W}_{n}=RD_{1}+\cdots +RD_{n}$, where $D_{i}=\frac{d}{%
dX_{i}}$. Thus $\mathbf{W}_{n}$ acts on the Grassmann algebra of K{{\"a}}hler
differential forms on $R$. We define three subalgebras of $\mathbf{W}_{n}.$

(1) The subalgebra $\mathbf{S}_{n}$, called the \emph{special algebra},
consists of those derivations annihilating the differential form $\nu
=dX_{1}\wedge \ldots \wedge dX_{n}.$ We can describe this set of derivations
as follows:
\[
\mathbf{S}_{n}=\left\{ \sum_{j=1}^{n}a_{j}D_{j}:\sum_{j=1}^{n}D_{j}\left(
a_{j}\right) =0\right\} .
\]

It is additively spanned by $D_{ij}(f)=D_{j}(f)D_{i}-D_{i}(f)D_{j}$, where $%
1\leq i,j\leq n$ and $f$ ranges over $R.$

(2) Suppose now that $n=2m$ is even then $\mathbf{H}_{2m},$ called the \emph{%
hamiltonian algebra}, consists of derivations annihilating $\omega
=\sum_{1\leq i\leq m}dX_{i}\wedge dX_{m+i}$. Define the following notation:
\begin{eqnarray*}
j^{\prime } &=&j+m\text{ and }\sigma (j)=1\text{ if }1\leq j\leq m \\
j^{\prime } &=&j-m\text{ and }\sigma (j)=-1\text{ if }m+1\leq j\leq 2m.
\end{eqnarray*}
Then we have the following description of the hamiltonian algebra:
\[
\mathbf{H}_{2m}=\left\{ \sum_{j=1}^{2m}a_{j}D_{j}:\sigma (j^{\prime
})D_{i}\left( a_{j^{\prime }}\right) =\sigma (i^{\prime })D_{j}\left(
a_{i^{\prime }}\right) ,1\leq i,j\leq 2m\right\} .
\]

Every element of $\mathbf{H}_{2m}$ can be represented as $D_{H}(f)$ for some
$f\in R$ where
\[
D_{H}(f)=\sum_{j=1}^{2m}\sigma (j)D_{j}(f)D_{j^{\prime }}.
\]

(3) Finally if $n=2m+1$ is odd then $\mathbf{K}_{2m+1},$ called the contact
algebra, consists of derivations $\partial $ such that $\partial \cdot
\alpha =f\cdot \alpha $ for some polynomial $f$ where
\[
\alpha =dX_{2m+1}+\sum_{1\leq i\leq m}X_{i}dX_{m+i}-X_{m+i}dX_{i}.
\]
For every derivation $D=\sum_{j=1}^{n}a_{j}\frac{d}{dX_{j}}$ define
\[
u(D)=\sum_{j=1}^{2m}\sigma (j)X_{j}\frac{da_{j^{\prime }}}{dX_{n}}-\frac{%
da_{n}}{dX_{n}}.
\]
Then $D$ is contained in $\mathbf{K}_{2m+1}$ if and only if for each $1\leq
i\leq 2m$%
\[
\sigma (i)a_{i^{\prime }}-\sum_{j=1}^{2m}\sigma (j)X_{j}\frac{da_{j^{\prime
}}}{dX_{i}}-\frac{da_{n}}{dX_{i}}=\sigma (i)X_{i^{\prime }}u(D).
\]

Every element of $\mathbf{K}_{2m+1}$ can be represented as $%
D_{K}(f)=\sum_{j=1}^{2m+1}f_{j}D_{j}$ where
\begin{eqnarray*}
f_{j} &=&X_{j}D_{2m+1}(f)+\sigma (j^{\prime })D_{j^{\prime }}(f)\text{ for }%
j\leq 2m \\
f_{2m+1} &=&2f-\sum_{j=1}^{2m}\sigma (j)X_{j}f_{j^{\prime }}.
\end{eqnarray*}

The Lie algebras $\mathbf{W}_{n},\mathbf{S}_{n},\mathbf{H}_{2m}$ and $%
\mathbf{K}_{2m+1}$ are called \emph{Lie algebras of Cartan type}. The usual
grading of the polynomial ring $k[X_{1},\ldots ,X_{n}]$ induces a grading of
$\mathbf{W}_{n},\mathbf{S}_{n}$ and $\mathbf{H}_{2m}.$ If $n=2m+1$ and $%
m\geq 1$ then there exists a unique grading of $k[X_{1},\ldots ,X_{2m+1}]$
such that $X_{1},\ldots ,X_{2m}$ are homogeneous of degree 1 and $X_{2m+1}$
is homogeneous of degree 2. The induced grading for $\mathbf{K}_{2m+1}$ is
called the natural grading of $\mathbf{K}_{2m+1}.$

\begin{theorem}
The Virasoro algebra $\mathrm{Vir}$ and the Lie
algebras of Cartan type $\mathbf{W}_{n},%
\mathbf{S}_{n},\mathbf{H}_{2m}$ and $\mathbf{K}_{2m+1}$ are well-covered.
\end{theorem}

\begin{proof}This depends on the fact that any derivation
\[
X_{1}^{m_{1}}\ldots X_{i}^{m_{i}+1}\ldots X_{n}^{m_{n}}\frac{d}{dX_{i}}
\]
with $m_{1}+\cdots +m_{n}=N$ can be realised as a linear multiple of the Lie
product
\[
\left( X_{1}^{r_{1}}\ldots X_{i}^{r_{i}+1}\ldots X_{n}^{r_{n}}\frac{d}{dX_{i}%
},X_{1}^{s_{1}}\ldots X_{i}^{s_{i}+1}\ldots X_{n}^{s_{n}}\frac{d}{dX_{i}}%
\right)
\]
where $r_{j}+s_{j}=m_{j}$ and $r_{i}\neq s_{i}$. For the Lie algebras
$\mathrm{Vir}$ and $\mathbf{W}_{n}$ the result follows, since
we can then realise the
derivation in the grading of weight $N$ as a Lie product of elements of
weight $i$ and $N-i$ for all values of $i<N/2.$

We shall have to work harder to get the corresponding result for the other
Cartan algebras $\mathbf{S}_{n},\mathbf{H}_{2m}$ or $%
\mathbf{K}_{2m+1}$,
since we need to realise elements of weight $N$ as sums of
commutators of elements of weight $i$ and $N-i$ which still stabilize the
form defining the corresponding Cartan algebra.

For $\mathbf{S}_{n}$ let us show why, for $\mathbf{X}^{\mathbf{c}%
}=X_{1}^{c_{1}}\ldots X_{n}^{c_{n}}$ where $\left\langle \mathbf{c}%
\right\rangle =c_{1}+\cdots +c_{n}=nN$, that $D_{ij}(\mathbf{X}^{\mathbf{c}})$
can be realised as a linear combination of commutators of elements of
weight $nN-t$ and $t$ for $N$ values of $t.$

We are going to use the following identity which can easily be checked.
Suppose that $%
b_{j}=0$ then
\begin{equation}
\left[ D_{ki}(\mathbf{X}^{\mathbf{a}}),D_{kj}(\mathbf{X}^{\mathbf{b}%
})\right] =a_{i}b_{k}D_{kj}(\mathbf{X}^{\mathbf{a}+\mathbf{b}-\mathbf{e}_{i}-%
\mathbf{e}_{k}})+a_{k}b_{k}D_{ij}(\mathbf{X}^{\mathbf{a}+\mathbf{b}-2\mathbf{%
e}_{k}})  \label{D_ij identity}
\end{equation}

Since $\left\langle \mathbf{c}\right\rangle =nN,$ we know that there is some
$c_{k}\geq N.$

Case 1. Suppose firstly that $k\neq i,j.$ Then choose $\mathbf{a}$ and $%
\mathbf{b}$ with $\mathbf{c}=\mathbf{a}+\mathbf{b}-2\mathbf{e}_{k},$ $%
b_{j}=0 $, $a_{i}=0$ and $a_{k}=c_{k}+1-l$ and $b_{k}=l+1$. For $l=0,\ldots
,c_{k}$, we get
\[
\left[ D_{ki}(\mathbf{X}^{\mathbf{a}}),D_{kj}(\mathbf{X}^{\mathbf{b}%
})\right] =a_{k}b_{k}D_{ij}(\mathbf{X}^{\mathbf{a}+\mathbf{b}-2\mathbf{e}%
_{k}})=\lambda D_{ij}(\mathbf{X}^{\mathbf{c}})
\]
with $\lambda \neq 0$ and $D_{ki}(\mathbf{X}^{\mathbf{a}})$ of weight $%
a_{1}+\cdots +a_{k-1}+a_{k+1}+\cdots +a_{n}+c_{k}+1-l$ and $D_{kj}(\mathbf{X}%
^{\mathbf{b}})$ of weight $b_{1}+\cdots +b_{k-1}+b_{k+1}+\cdots +b_{n}+l+1.$

Case 2. If $k=i,$ then choose some $s\in \left\{ 1,\ldots ,n\right\}
\backslash \left\{ i,j\right\} $ and define $\mathbf{a}$ and $\mathbf{b}$
with $\mathbf{c}=\mathbf{a}+\mathbf{b}-2\mathbf{e}_{s},$ $b_{j}=0$, $a_{s}=1$
and $a_{i}=c_{i}+1-l$ and $b_{i}=l+1$. Form
\[
\left[ D_{si}(\mathbf{X}^{\mathbf{a}}),D_{sj}(\mathbf{X}^{\mathbf{b}%
})\right] =a_{i}b_{s}D_{sj}(\mathbf{X}^{\mathbf{a}+\mathbf{b}-\mathbf{e}_{i}-%
\mathbf{e}_{s}})+a_{s}b_{s}D_{ij}(\mathbf{X}^{\mathbf{a}+\mathbf{b}-2\mathbf{%
e}_{s}})
\]
where note that $a_{s}\neq 0\neq b_{s}.$ To deal with $a_{i}b_{s}D_{sj}(%
\mathbf{X}^{\mathbf{a}+\mathbf{b}-\mathbf{e}_{i}-\mathbf{e}_{s}})$ we can
use the analysis of case 1 since the weight is concentrated still at $c_{i}$
but now we are considering a derivation $D_{sj}$ with $i\neq s,j.$ Hence
choose $\mathbf{a}^{\prime }=\mathbf{a+e}_{i}-\mathbf{e}_{s}$ then $\mathbf{a%
}+\mathbf{b}-\mathbf{e}_{i}-\mathbf{e}_{s}=\mathbf{a}^{\prime }+\mathbf{b}-2%
\mathbf{e}_{i}$, and $a_{s}^{\prime }=0$, $a_{i}^{\prime }=c_{i}+2-l$ and $%
b_{i}=l+1$ where $l=0,\ldots ,c_{i}.$ (Note that the weight at $X_{i}$ of $%
D_{sj}(\mathbf{X}^{\mathbf{a}+\mathbf{b}-\mathbf{e}_{i}-\mathbf{e}_{k}})$ is
now $c_{i}+1.)$ Hence we have
\[
\left[ D_{si}(\mathbf{X}^{\mathbf{a}}),D_{sj}(\mathbf{X}^{\mathbf{b}%
})\right] -(a_{i}b_{s})(a_{i}^{\prime }b_{i}^{\prime })^{-1}\left[ D_{is}(%
\mathbf{X}^{\mathbf{a}^{\prime }}),D_{ij}(\mathbf{X}^{\mathbf{b}})\right]
=\lambda D_{ij}(\mathbf{X}^{\mathbf{c}})
\]
with $\lambda \neq 0,$ $D_{si}(\mathbf{X}^{\mathbf{a}})$ and $D_{is}(\mathbf{%
X}^{\mathbf{a}^{\prime }})$ have weight $a_{1}+\cdots
+a_{i-1}+a_{i+1}+\cdots +a_{n}+c_{i}+1-l$ and $D_{sj}(\mathbf{X}^{\mathbf{b}
})$ and $D_{ij}(\mathbf{X}^{\mathbf{b}})$ have weight $b_{1}+\cdots
+b_{i-1}+b_{i+1}+\cdots +b_{n}+l+1.$ This proves then that $\mathbf{S}_{n}$
is well-covered.

For the Hamiltonian algebra, we can use the following identity (see Chapter
4 Lemma 4.3 (1) of \cite{Strade and Farnsteiner})
\begin{eqnarray*}
\left[ D_{H}(\mathbf{X}^{\mathbf{a}}),D_{H}(X_{i}^{b_{i}})\right]
&=&D_{H}\left( D_{H}(\mathbf{X}^{\mathbf{a}})(X_{i}^{b_{i}})\right) \\
&=&D_{H}\left( \sum_{j=1}^{2m}\sigma (j^{\prime })D_{j^{\prime
}}(f)D_{j}(X_{i}^{b_{i}})\right) \\
&=&D_{H}\left( \sigma (i^{\prime })a_{i^{\prime }}\mathbf{X}^{\mathbf{a}-%
\mathbf{e}_{i^{\prime }}}b_{i}X_{i}^{b_{i}-1}\right) .
\end{eqnarray*}
Consider now realising $D_{H}(\mathbf{X}^{\mathbf{c}})$ of weight $Nn.$ So
there exists some $i$ such that $c_{i}\geq N.$ We then put $b_{i}=l+1$ and $%
\mathbf{a}=\mathbf{c}+\mathbf{e}_{i^{\prime }}-b_{i}+1$ with $l=0,\ldots
,c_{i}$ and hence can express $D_{H}(\mathbf{X}^{\mathbf{c}})$ as a
commutator of elements of the Hamiltonian algebra of weight $Nn-l$ and $l$
for $l=0,\ldots ,c_{i}$. Hence $\mathbf{H}_{2m}$ is well-covered.

For the contact algebra $\mathbf{K}_{n}=\mathbf{K}_{2m+1}$, we shall use the
following identity (see Chapter 4 Proposition 5.2 of \cite{Strade and
Farnsteiner}):
\[
\left[ D_{K}(f),D_{K}(g)\right] =D_{K}\left(
D_{K}(f)(g)-2gD_{2m+1}(f)\right) .
\]

This algebra is slightly trickier than the previous algebras. We shall need
to consider realising $D_{K}(\mathbf{X}^{\mathbf{c}})$ of weight $4Nn(n-1)$
as a linear combination of commutators of weight $4Nn(n-1)-l$ and $l$ for $N$
different values of $l.$ Recall that in this algebra $X_{n}$ has weight two
whilst the other variables have weight one. There exists some $i$ with $c_{i}%
\geq 2N(n-1).$

Case 1. Suppose that $i\neq n.$ We prove by induction on $c_{n}$ that we can
express $D_{K}(\mathbf{X}^{\mathbf{c}})$ as a linear combination of
commutators of weight $4Nn(n-1)-l$ and $l$ for $l=0,\ldots ,N.$ Put $f=%
\mathbf{X}^{\mathbf{a}}$ and $g=X_{i}^{b_{i}}$ where $\mathbf{a}=\mathbf{c+e}%
_{i^{\prime }}-l\mathbf{e}_{i}$ and $b_{i}=l+1.$ Then
\begin{equation}
\left[ D_{K}(\mathbf{X}^{\mathbf{a}}),D_{K}(X_{i}^{b_{i}})\right]
=D_{K}\left( X_{i}a_{n}\mathbf{X}^{\mathbf{a}-\mathbf{e}%
_{n}}b_{i}X_{i}^{b_{i}-1}\right) +D_{K}\left( \sigma (i^{\prime
})a_{i^{\prime }}\mathbf{X}^{\mathbf{a}-\mathbf{e}_{i^{\prime
}}}b_{i}X_{i}^{b_{i}-1}\right)   \label{D_K identity}
\end{equation}
Suppose that $c_{n}=a_{n}=0.$ Then we are done. We then suppose by induction
that we have proved the claim for $c_{n}-1.$ Then the identity (\ref{D_K
identity}) shows how to express $D_{K}(\mathbf{X}^{\mathbf{c}})$ as a linear
combination of commutators as desired since the first expression on the
right hand side of (\ref{D_K identity}) can be dealt with by the inductive
hypothesis.

Case 2. Suppose now that $i=n.$ We show how to shift the weight from $X_{n}$
to the other variables such that eventually we are in a position to invoke
case 1.

We use the following identity
\[
\left[ D_{K}(X_{n}^{b_{n}}),D_{K}(\mathbf{X}^{\mathbf{a}})\right]
=D_{K}\left( \left\langle X_{n}^{b_{n}},\mathbf{X}^{\mathbf{a}}\right\rangle
\right)
\]
where
\begin{equation*}
\begin{split}
\left\langle X_{n}^{b_{n}},\mathbf{X}^{\mathbf{a}}\right\rangle
&=D_{K}(X_{n}^{b_{n}})(\mathbf{X}^{\mathbf{a}})-2\mathbf{X}^{\mathbf{a}%
}D_{n}(X^{b_{n}}) \\
&=\left( \left( \left( \sum_{j=1}^{2m}a_{j}\right) -2\right)
b_{n}+2a_{n}\right) \mathbf{X}^{\mathbf{a}}X_{n}^{b_{n}-1}\\
&-\sum_{j=1}^{2m}%
\sigma (j)b_{n}a_{n}\mathbf{X}^{\mathbf{a}+(b_{n}-2)\mathbf{e}_{n}+\mathbf{e}%
_{j}+\mathbf{e}_{j^{\prime }}}.\\
\end{split}
\end{equation*}

Note that we shall let $b_{n}=l+1$ for $l=0,\ldots ,N$ and $a_{n}=c_{n}-l.$
Since $c_{n}\geq 2N(n-1)$ this will mean that  $\left( \left( \left(
\sum_{j=1}^{2m}a_{j}\right) -2\right) b_{n}+2a_{n}\right) $ is always
non-zero. The choice of $c_{n}\geq 2N(n-1)$ means that we can keep on
applying the above identity to realise the expressions $\sigma (j)b_{n}a_{n}%
\mathbf{X}^{\mathbf{a}+(b_{n}-2)\mathbf{e}_{n}+\mathbf{e}_{j}+\mathbf{e}%
_{j^{\prime }}}$ as commutators $\left[ D_{K}(X_{n}^{b_{n}^{\prime }}),D_{K}(%
\mathbf{X}^{\mathbf{a}^{\prime }})\right] $ with an error term which has the
values of $c_{i}$ increasing for $i\neq n$ and $c_{n}$ decreasing whilst
still ensuring that $c_{n}>N$ so that we can choose $b_{n}^{\prime }=l+1$
for $l=0,\ldots ,N.$ Eventually the error terms will have some $c_{i}$ with $%
c_{i}>N$ and we can apply case 1 to finish the realisation.
\end{proof}

\begin{corollary}
If $L$ is a simple $\Bbb{Z}$-graded Lie algebra over an algebraically closed
field of characteristic zero of finite growth then $\mathcal{X}_{0}^{*}(L)$
and $\mathcal{X}^{*}(L_{0})$ for $*\in \left\{ \leq ,\triangleleft \right\} $
are stable classes of subalgebras.
\end{corollary}

The proof above depends on the classification proved by Mathieu
\cite{Mathieu}
and a case
by case analysis of each class of Lie algebra in the classification. It may
be possible that there is a more direct argument based on the simplicity of
the $\Bbb{Z}$-graded Lie algebra. Note that since $L$ is simple $\mathcal{X}
_{0}^{\triangleleft }(L)$ is actually empty or consists of $L$ if $L_{0}$
has finite codimension in $L$.

As we have pointed out, an $\Bbb{N}$-filtered Lie algebra has many ideals so
we cannot expect it to be simple. However in this context, the concept of
being simple is replaced by that of being \emph{just infinite} - that is an
infinite dimensional Lie algebra all of whose proper quotients are finite
dimensional. It is conjectured by Shalev and Zelmanov (see \cite{Shalev-New
Horizons} section 6.5) that the just infinite $\Bbb{N}$-filtered Lie
algebras over an algebraically closed field of characteristic zero for which
$L_{0}=L_{1}$ and $\dim (L_{i}/L_{i+1})$ are uniformly bounded are in the
following list:

(1) $L$ is soluble;

(2) The completion of $L$ is commensurable with the positive part of a loop
algebra;

(3) The completion of $L$ is commensurable to the completion of $\partial
^{+}.$

Note that only (2) and (3) will have stable subalgebras.

\subsection{Non-stability of free graded Lie algebras}

Stability for infinite Lie algebras is the analogue of determining when a
group has only a finite number of subgroups of each given index which can
then be counted using a standard Dirichlet series. For groups the condition
that the group be finitely generated (either as an abstract group or a
topological group) was sufficient to ensure this condition. Note that being
finitely generated will not suffice in the context of infinite dimensional
Lie algebras as the following Theorem indicates:

\begin{theorem}
Let $L=\bigoplus_{i\geq 1}\gamma _{i}(L)/\gamma _{i+1}(L)$ be the free two
generated infinite dimensional graded Lie algebra where $\gamma _{i}(L)$ is
the lower central series. Then for each $i$ there exist closed subalgebras of
codimension 2 not containing $\gamma _{i}(L).$
\end{theorem}

\begin{proof}Let $x$ and $y$ be the free generators. If $w$ is any Lie word
in $x$ and $y$ we define the length $l(w)$ to be the number of terms in $w.$
We take a Hall set $H=\left\{ w_{i}:i\geq 1\right\} $ as a basis for $L$
(see \cite{Bourbaki-Lie algebra} II.2.10) which is defined as follows:

(1) if $w_{i}\in H$ and $w_{j}\in H$ and $l(w_{i})<l(w_{j})$ then $i<j;$

(2) $w_{1}=x,w_{2}=y$ and $w_{3}=(xy);$

(3) an element $w$ of length $\geq 3$ belongs to $H$ if and only if it is of
the form $(a(bc))$ with $a,b,c$ in $H$, $(bc)\in H$ and $b\leq a<bc$ and $
b<c $ where the ordering is defined by the ordering on the index set.

For example the construction provided by Proposition 11 of \cite
{Bourbaki-Lie algebra} II.2.10 starts with the following basis for each
layer $\gamma_{i}(L)/\gamma_{i+1}(L)$:
\[
\begin{array}{llll}
\gamma_{1}(L)/\gamma_{2}(L) & w_{1}=x & w_{2}=y &  \\
\gamma_{2}(L)/\gamma_{3}(L) & w_{3}=(xy) &  &  \\
\gamma_{3}(L)/\gamma_{4}(L) & w_{4}=(x(xy)) & w_{5}=(y(xy)) &  \\
\gamma_{4}(L)/\gamma_{5}(L) & w_{6}=(x(x(xy))) & w_{7}=(y(x(xy))) &
w_{8}=(y(y(xy))) \\
\gamma_{5}(L)/\gamma_{6}(L) & w_{9}=(x(x(x(xy)))) & w_{10}=(y(x(x(xy)))) &
w_{11}=(y(y(x(xy)))) \\
& w_{12}=(y(y(y(xy)))) & w_{13}=((xy)(x(xy))) & w_{14}=((xy)(y(xy)))
\end{array}
\]

Let $n_{i}-1$ be the dimension of $\gamma_{1}(L)/\gamma_{i}(L).$ Then we
claim that the additive subspace of $L$ generated by the following basis
elements is actually a subalgebra:
\[
H_{i}=\bigoplus_{l\neq 1,n_{i}}{k}w_{l}.
\]
So this is a vector subspace of codimension 2 which skips the basis elements
$w_{1}=x$ and $w_{n_{i}}=(x(x\ldots(xy)\ldots))$ where $x$ appears $i-1$
times.

The Hall basis has the property that each word $w_{l}$ has a unique
decreasing factorization $w_{l}=(w_{j}w_{k})$ with $j<k<l$ (see Corollary
4.7 of \cite{Reutenauer}).

We prove by induction on the length of $w_{j}$ for $1<j<n_{i}$ and $k\neq
1,n_{i}$ that $(w_{j}w_{k})=\sum_{l\neq 1,n_{i}}a_{l}w_{l}.$ If $w_{j}$ has
length $1$ then $w_{j}=y$ and by condition (3) $(yw_{k})$ is in $H.$ The
property of unique factorization implies it is not the element $%
w_{n_{i}}=(xw_{n_{i-1}}).$ Suppose we have proved that $(w_{j}w_{k})=\sum_{l%
\neq 1,n_{i}}a_{l}w_{l}$ for all words $w_{j}$ of length less that $m$ and
take $w_{j}$ of length $m>1.$ There is a unique decreasing factorization $
w_{j}=(w_{j_{1}}w_{j_{2}})$ with $l(w_{j_{i}})<l(w_{j}).$ Using the Jacobi
identity we can rewrite
\begin{eqnarray*}
(w_{j}w_{k}) &=&((w_{j_{1}}w_{j_{2}})w_{k}) \\
&=&(w_{j_{1}}(w_{j_{2}}w_{k}))-(w_{j_{2}}(w_{j_{1}}w_{k})).
\end{eqnarray*}
Now we can use our induction hypothesis applied to $w_{j_{1}}$ and $%
w_{j_{2}}.$ \end{proof} Obviously the same argument applies to the
$d$-generated free Lie algebra.

So the existence of a motivic zeta function for an infinite dimensional Lie
algebra does not appear to be so straightforward. The condition of being
well-covered will suffice, but this is a special sort of property which
depends on the internal workings of a Lie algebra. It would be nice to have
a more transparent Lie theoretic condition which will ensure the stability
of the class of subalgebras. We therefore raise the following:

\begin{problem}Determine a criterion for the class of all subalgebras of
finite codimension to be stable in an ${\Bbb N}$-filtered or ${\Bbb Z}$%
-filtered Lie algebra.
\end{problem}

>From now on we shall focus on the situation of a finite dimensional $k[[t]]$%
-Lie algebra and ${\cal X}$ a class of $k[[t]]$-subalgebras.

\section{Motivic integration and rationality of Poincar{\'e} series}

In this section, we review material from \cite{Denef-Loeser-Arcs} and \cite
{MK} which will be used in the present paper.

\subsection{Scheme of arcs}

For $X$ a variety over $k$, we will denote by $\cL (X)$ the {\em scheme of
germs of arcs on $X$}. It is a scheme over $k$ and for any field extension $%
k \subset K$ there is a natural bijection
\[
\cL (X) (K) \simeq {\rm Mor}_{k-{\rm schemes}} (\Spec K [[t]], X)
\]
between the set of $K$-rational points of $\cL (X)$ and the set of germs of
arcs with coefficients in $K$ on $X$. We will call $K$-rational points of $%
\cL (X)$, for $K$ a field extension of $k$, arcs on $X$, and $\varphi (0)$
will be called the origin of the arc $\varphi$. More precisely the scheme $%
\cL (X)$ is defined as the projective limit
\[
\cL (X) := \varprojlim \cL_{n} (X)
\]
in the category of $k$-schemes of the schemes $\cL_{n}(X)$ representing the
functor
\[
R \mapsto {\rm Mor}_{k-{\rm schemes}} (\Spec R [t] / t^{n+1} R[t], X)
\]
defined on the category of $k$-algebras. The existence of $\cL_{n}(X)$ is
well known (cf. \cite{Denef-Loeser-Arcs}) and the projective limit exists
since the transition morphisms are affine. We shall denote by $\pi_{n}$ the
canonical morphism, corresponding to truncation of arcs,
\[
\pi_{n} : \cL (X) \longrightarrow \cL_{n} (X).
\]
The schemes $\cL (X)$ and $\cL_{n} (X)$ will always be considered with their
reduced structure.

\subsection{Semi-algebraic geometry}

>From now on we will denote by $\bar{k}$ a fixed algebraic closure of $k$,
and by $\bar{k}((t))$ the fraction field of $\bar{k}[[t]]$, where $t$ is one
variable. Let $x_{1},\ldots ,x_{m}$ be variables running over $\bar{k}((t))$
and let $\ell _{1},\ldots ,\ell _{r}$ be variables running over ${\Bbb Z}$.
A {\it semi-algebraic} (resp. $k [[t]]$-{\it semi-algebraic}) condition $\theta
(x_{1},\ldots ,x_{m};\ell _{1},\ldots ,\ell _{r})$ is a finite boolean
combination of conditions of the form
\begin{align*}
\text{(1)} \quad & \ord_{t}f_{1}(x_{1},\ldots ,x_{m})\geq \ord_{t}f_{2}(x_{1},%
\ldots ,x_{m})+L(\ell _{1},\ldots ,\ell _{r}) \\
\text{(2)} \quad& \ord_{t}f_{1}(x_{1},\ldots ,x_{m})\equiv L(\ell _{1},\ldots
,\ell _{r})\pmod d \\
\text{(3)} \quad & h(\overline{ac}(f_{1}(x_{1},\ldots ,x_{m})),\ldots
,\overline{ac}%
(f_{m^{\prime }}(x_{1},\ldots ,x_{m})))=0,
\end{align*}
where $f_{i}$ are polynomials with coefficients in $k$ (resp. $f_{i}$ are
polynomials with coefficients in $k[[t]]$), $h$ is a polynomial with
coefficients in $k$, $L$ is a polynomial of degree $\leq 1$ over ${\Bbb Z}$,
$d\in {\Bbb N}$, and $\overline{ac}(x)$ is the coefficient of lowest degree
of $x$ in $\bar{k}((t))$ if $x \not= 0$, and is equal to 0 otherwise. Here
we use the convention that $\infty +\ell =\infty $ and $\infty \equiv \ell \;%
{\rm mod}\,d$, for all $\ell \in {\Bbb Z}$. In particular the condition $%
f(x_{1},\ldots ,x_{m})=0$ is a semi-algebraic condition
(resp. a $k[[t]]$-semi-algebraic condition), for $f$ a
polynomial over $k$ (resp. over $k[[t]]$). A subset of
$\bar{k}((t))^{m}\times {\Bbb Z}^{r}$
defined by a semi-algebraic (resp. $k [[t]]$-semi-algebraic) condition is called
{\em semi-algebraic} (resp. {\em $k [[t]]$-semi-algebraic}).
One defines similarly
semi-algebraic and $k [[t]]$-semi-algebraic subsets of $K((t))^{m}\times
{\Bbb Z}%
^{r}$ for $K$ an algebraically closed field containing $\bar{k}$.

A function $\alpha : \bar k ((t))^m \times {\Bbb Z}^n \rightarrow {\Bbb Z}$
is called {\it simple} (resp. $k [[t]]$-{\it simple}) if its graph is
semi-algebraic (resp. $k [[t]]$-semi-algebraic).

Let $X$ be an algebraic variety over $k$. For $x\in \cL (X)$, we denote by $%
k_{x}$ the residue field of $x$ on $\cL (X)$, and by $\tilde{x}$ the
corresponding rational point $\tilde{x}\in \cL (X)(k_{x})=X(k_{x}[[t]])$.
When there is no danger of confusion we will often write $x$ instead of $%
\tilde{x}$. A {\it semi-algebraic family of semi-algebraic subsets}
(resp.
$k [[t]]$-{\it semi-algebraic family of $k [[t]]$-semi-algebraic subsets})
(for $n=0$ a
semi-algebraic subset (resp. $k [[t]]$-semi-algebraic subset)) $A_{\ell }$,
$\ell
\in {\Bbb N}^{n}$, of $\cL (X)$ is a family of subsets $A_{\ell }$ of $\cL %
(X)$ such that there exists a covering of $X$ by affine Zariski open sets $U$
with
\[
A_{\ell }\cap \cL (U)=\Bigl\{x\in \cL (U) :
\theta (h_{1}(\tilde{x}),\ldots ,h_{m}(\tilde{x});\ell )\Bigr\},
\]
where $h_{1},\ldots ,h_{m}$ are regular functions on $U$ and $\theta $ is a
semi-algebraic condition (resp. $k [[t]]$-semi-algebraic condition). Here the $%
h_{i}$'s and $\theta $ may depend on $U$ and $h_{i}(\tilde{x})$ belongs to $
k_{x}[[t]]$.

Let $A$ be a semi-algebraic subset (resp. $k [[t]]$-semi-algebraic subset) of $%
\cL
(X)$. A function $\alpha : A \times {\Bbb Z}^n \rightarrow {\Bbb Z} \cup
\{\infty \}$ is called {\it simple} (resp. $k [[t]]$-{\it simple}) if the family
of subsets $\{x \in \cL (X) : \alpha (x, \ell_1, \ldots, \ell_n) =
\ell_{n + 1}\}$, $(\ell_1, \ldots, \ell_{n + 1}) \in {\Bbb N}^{n + 1}$, is a
semi-algebraic family of semi-algebraic subsets (resp. a $k [[t]]$-{\it
semi-algebraic family of $k [[t]]$-semi-algebraic subsets}) of $\cL (X)$.

An important fact is that if $A$ is a
$k [[t]]$-semi-algebraic subset of $\cL (X)$,
then $\pi_n (A)$ is a constructible subset of $\cL_n (X)$ (cf. Proposition
1.7 of \cite{MK}).

\begin{remark}Motivic integration
is developed in \cite{Denef-Loeser-Arcs}
for semi-algebraic subsets of $\cL (X)$, when $X$ is an algebraic variety
over $k$,
and simple functions. In the paper
\cite{MK}, this is
 extended to what is called there $t$-semi-algebraic subsets of $\cL (X)$
and $t$-simple functions, which are defined similarly as
$k [[t]]$-semi-algebraic subsets of $\cL (X)$
and $k [[t]]$-simple functions, but replacing $k [[t]]$ by $k [t]$.
In fact,
as mentioned in Remark 1.18 of \cite{MK}, all results in
\S\kern .15em 1 of \cite{MK} (before 1.17) may be extended, with
similar proofs, to cover the case of
$k [[t]]$-semi-algebraic subsets
and $k [[t]]$-simple functions. Hence when we shall quote
a result from \cite{MK}, we shall use its extension to
$k [[t]]$-semi-algebraic subsets
and $k [[t]]$-simple functions without further comment.
\end{remark}

\subsection{Motivic integration}

>From now on we assume $X$ is a smooth algebraic variety over $k$ of pure
dimension $d$. Let $A$ be a $k [[t]]$-semi-algebraic subset $\cL (X)$. We
say $A$
is stable at level $n$ if $A=\pi _{n}^{-1}\pi _{n}(A)$. Remark that if $A$
is stable at level $n$, then $A$ is stable at level $m$, for any $m\geq n$.
We say $A$ is stable if it is stable at some level. Denote by $\BB^{t}$ the
set of all $k [[t]]$-semi-algebraic subsets of $\cL (X)$, and by
$\BB_{0}^{t}$ the
set of all $A$ in $\BB^{t}$ which are stable. Clearly there is a unique
additive measure
\[
\tilde{\mu}:\BB_{0}^{t}\longrightarrow \cM_{{\rm loc}}
\]
satisfying
\[
\tilde{\mu}(A)=[\pi _{n}(A)]\,\LL^{-(n+1)d}\,,
\]
when $A$ is stable at level $n$. In fact, the condition of being $k [[t]]$%
-semi-algebraic is superfluous here. By the same formula one may define $%
\tilde{\mu}(A)$ for $A$ cylindrical at level $n$, {\it i.e.} subsets of $\cL
(X)$ of the form $A=\pi _{n}^{-1}(C)$ with $C$ constructible. One says $A$
is cylindrical if it is cylindrical at some level.

Let $A$ be in $\BB_{0}^{t}$ and let $\alpha ~:A\rightarrow {\Bbb Z}$
be a $k [[t]]$%
-simple function (or, more generally, assume $A$ and the fibers of $\alpha $
are cylindrical). By Lemma 2.4 of \cite{Denef-Loeser-Arcs} and Lemma A.3 of
\cite{MK}, $|\alpha |$ is bounded, and we can define
\begin{equation}\label{int}
\int_{A}\LL^{-\alpha }d\tilde{\mu}~:=\sum \LL^{-n}\tilde{\mu}(\alpha
^{-1}(n)),
\end{equation}
the sum on the right hand side being finite.

\subsection{Completion}

\label{completion}We now explain how one extends $\tilde \mu$ to non
stable
$k [[t]]$-semi-algebraic subsets by using a completion of $\cM_{{\rm
loc}}$. This
is indeed similar to the use of real numbers for defining $p$-adic
integrals. The material here will only be used in section \ref
{Section:Motivic cone integrals}. So we denote by $\widehat \cM$ the
completion of ${\cal M}_{{\rm loc}}$ with respect to the filtration $F^{m}%
{\cal M}_{{\rm loc}}$ where $F^{m}{\cal M}_{{\rm loc}}$ is the subgroup
generated by $\left\{ \left[ S\right] {\bf L}^{-i}:i-\dim S\geq m\right\}$.
We will also denote by $F^{\cdot}$ the filtration induced on $\widehat{\cM}$.

In \cite{Denef-Loeser-Arcs} and \cite{MK} the following is shown\footnote{%
in fact loc. cit. also covers the case of singular varieties}: There exists
a unique map $\mu : \BB^{t} \rightarrow \widehat \cM$ satisfying the
following three properties.

\begin{enumerate}
\item[(1)]  If $A\in \BB^{t}$ is stable at level $n$, then $\mu (A)=[\pi %
_{n}(A)]\LL^{-(n+1)d}$.

\item[(2)]  If $A\in \BB^{t}$ is contained in $\cL (S)$ with $S$ a reduced
closed subscheme of $X\otimes k[[t]]$ with ${\rm dim}_{k[[t]]}\,S<{\rm dim}\,X$,
then $\mu (A)=0$.

\item[(3)]  Let $A_{i}$ be in $\BB^{t}$ for each $i$ in ${\Bbb N}$. Assume
that the $A_{i}$'s are mutually disjoint and that $A:=\bigcup_{i\in {\Bbb N}%
}A_{i}$ is $k [[t]]$-semi-algebraic. Then $\sum_{i\in {\Bbb N}}\mu (A_{i})$
converges in $\widehat{\cM}$ to $\mu (A)$.
\end{enumerate}

Moreover we have:

\begin{enumerate}
\item[(4)]  If $A$ and $B$ are in $\BB^{t}$, $A\subset B$ and if $\mu (B)\in
F^{m}\widehat{\cM}$, then $\mu (A)\in F^{m}\widehat{\cM}$.
\end{enumerate}

This unique map $\mu$ is called the {\it motivic volume} on $\cL (X)$ and is
denoted by $\mu_{\cL (X)}$ or $\mu$. For $A$ in $\BB^{t}$ and $\alpha : A
\rightarrow {\Bbb Z}\cup \{\infty\}$ a $k [[t]]$-simple function, one
defines the
motivic integral
\begin{equation}\label{int2}
\int_{A} \LL^{- \alpha} d \mu := \sum_{n \in {\Bbb Z} } \mu (A \cap
\alpha^{-1} (n)) \, \LL^{- n}
\end{equation}
in $\widehat \cM$, whenever the right hand side converges in $\widehat \cM$,
in which case we say that $\LL^{- \alpha}$ is integrable on $A$. If the
function $\alpha$ is bounded from below, then $\LL^{- \alpha}$ is integrable
on $A$, because of (4).

\subsection{Rationality}

The following results are proved in section 5 of \cite{Denef-Loeser-Arcs}
for semi-algebraic families of semi-algebraic subsets and simple functions.
For simplicity results are stated here only for smooth varieties.

\begin{theorem}
\label{5.1prime}Let $X$ be a smooth algebraic variety over $k$ of pure
dimension $d$. Let $A_{n}$, $n\in {\Bbb Z}^{r}$, be a semi-algebraic family
of semi-algebraic subsets of $\cL (X)$ and let $\alpha :\cL (X)\times {\Bbb Z%
}^{r}\rightarrow {\Bbb N}$ be a simple function. Assume that $A_{n}$ and the
fibers of $\alpha (\_,n):A_{n}\rightarrow {\Bbb N}$ are stable, for every $%
n\in {\Bbb N}^{r}$. Then the power series
\[
\sum_{n\in {\Bbb N}^{r}}T^{n}\,\int_{A_{n}}\LL^{-\alpha (\_,n)}d\tilde{\mu}
\]
in the variable $T=(T_{1},\ldots ,T_{r})$ belongs to the subring of $\cM_{%
{\rm loc}}[[T]]$ generated by $\cM_{{\rm loc}}[T]$ and the series $(1-\LL%
^{-a}T^{b})^{-1}$, with $a\in {\Bbb N}$ and $b\in {\Bbb N}^{r}\setminus \{0\}
$.
\end{theorem}

\begin{theorem}
\label{5.1}Let $X$ be a smooth algebraic variety over $k$ of pure dimension $%
d$. Let $A_{n}$, $n\in {\Bbb Z}^{r}$, be a semi-algebraic family of
semi-algebraic subsets of $\cL (X)$ and let $\alpha :\cL (X)\times {\Bbb Z}
^{r}\rightarrow {\Bbb N}$ be a simple function. Then the power series
\[
\sum_{n\in {\Bbb N}^{r}}T^{n}\,\int_{A_{n}}\LL^{-\alpha (\_,n)}d\mu
\]
in the variable $T=(T_{1},\ldots ,T_{r})$ belongs to the subring of $%
\widehat{\cM }[[T]]$ generated by the image in $\widehat{\cM }[[T]]$ of $\cM
_{{\rm loc}}[T]$, $(\LL^{i}-1)^{-1}$ and $(1-\LL^{-a}T^{b})^{-1}$, with $%
i\in {\Bbb N}\setminus \{0\}$, $a\in {\Bbb N}$, $b\in {\Bbb N}^{r}\setminus
\{0\}$.
\end{theorem}

\begin{corollary}
{}For any semi-algebraic subset $A$ of $\cL (X)$, the measure $\mu (A)$ is
in $\overline{\cM}_{{\rm loc}}[((\LL^{i}-1)^{-1})_{i\geq 1}]$, where $%
\overline{\cM}_{{\rm loc}}$ is the image of $\cM_{{\rm loc}}$ in $\widehat{%
\cM}$ \hfill $\qed$
\end{corollary}

\begin{remark}
\label{min}By replacing $T_{i}$ by $T_{i}\LL^{m_{i}}$ in Theorems \ref
{5.1prime} and \ref{5.1}, one sees
that the condition ``$\alpha $ takes
values in ${\Bbb N}$'' may be replaced by the condition ``$\alpha $ is
bounded from below by a linear function of the ${\Bbb Z}^{r}$-variable''.
\end{remark}

\begin{remark}
\label{bad}It seems quite likely that Theorems \ref{5.1prime} and \ref{5.1}
remain true if semi-algebraic is replaced everywhere by $k [[t]]$-semi-algebraic
and simple by $k [[t]]$-simple. However it does not seem that the proofs may be
adapted directly to that more general situation, since there might be some
``bad reduction at $t=0$''.
\end{remark}

\section{Motivic integration on the infinite Grassmannian\label%
{section:Grassmannian}}

All the material in \ref{4.1} and \ref{4.2} is contained in \cite{BL}, or is
directly adapted, replacing ${\rm SL}$ by ${\rm GL}$, from statements in
\cite{BL}.

\subsection{}\label{4.1} We work over a field $k$ of characteristic $0$. We
shall
consider the infinite Grassmannian as a functor on the category of
$k$-algebras. (A $k$-algebra will always be assumed to be associative,
commutative and unitary.) A natural framework is that of $k$-spaces
and $k$-groups in the sense of \cite{BL}. By definition, a $k$-space (resp.
a $k$-group) is a functor from the category of $k$-algebras to the category
of sets (resp. of groups) which is a sheaf for the faithfully flat topology
(see \cite{BL} for a definition). The category of schemes over $k$ can be
viewed as a full subcategory of the category of $k$-spaces. Schemes will
always be assumed to be quasi-compact and quasi-separated. An important
feature is that direct limits exist in the category of $k$-spaces, so we can
say a $k$-space (resp. a $k$-group) is an {\it ind-scheme} (resp. an {\it
ind-group}) if it is the direct limit of a directed system of schemes (resp.
of group schemes).

We fix a positive integer $d$. Let $t$ be an indeterminate. We consider the $%
k$-groups ${\bf GL}_{d}(k[[t]])$ and ${\bf GL}_{d}\bigl(k((t))\bigr)$
respectively defined by $R\mapsto {\rm GL}_{d}(R[[t]])$ and $R\mapsto {\rm GL%
}_{d}\bigl(R((t))\bigr)$.

For $n\ge 0$, we denote by $G_{(n)}(R)$ the set of matrices $A(t)$ in ${\rm %
GL}_{d}\bigl(R((t))\bigr)$ such that both $A(t)$ and $A(t)^{-1}$ have a pole
of order $\leq n$. This defines a subfunctor $G_{(n)}$ of ${\bf GL}
_{d}(k((t)))$. One can show (cf. \cite{BL}), that the $k$-group ${\bf GL}
_{d}(k[[t]])$ is an affine group scheme and that the $k$-group ${\bf GL}
_{d}(k((t)))$ is an ind-group, being the direct limit of the sequence of
schemes $%
(G_{(n)})$, $n\geq 0$.

\subsection{}\label{4.2}The infinite Grassmannian $\gr$ may be defined as
the quotient
\[
{\bf GL}_{d}(k((t)))/{\bf GL}_{d}(k[[t]])
\]
in the category of $k$-spaces. It is the sheaf, for the faithfully flat
topology, associated to the presheaf $R\mapsto {\rm GL}_{d}(R((t)))/{\rm GL}
_{d}(R[[t]])$. Now we have to explain why $\gr$ is indeed the infinite
Grassmannian. For any $k$-algebra $R$, consider the set $W(R)$ of $R[[t]]$%
-submodules $L$ of $R((t))^{d}$ such that, for some $n\geq 0$, $%
t^{n}R[[t]]^{d}\subset L\subset t^{-n}R[[t]]^{d}$ and $L/t^{n}R[[t]]^{d}$ is
a projective $R$-module. By Proposition 2.3 of \cite{BL}, the $k$-space $\gr$
is isomorphic to the functor $R\mapsto W(R)$. Under this isomorphism, the
group scheme action of ${\bf GL}_{d}(k((t)))$ on $\gr$ corresponds to the
natural action of ${\rm GL}_{d}(R((t)))$ on $R[[t]]$-submodules of $%
R((t))^{d}$. Denote by ${\cal Gr}_{(n)}$ the image of $G_{(n)}$ in ${\cal Gr}
$. Under the preceding isomorphism ${\cal Gr}_{(n)}$ may be identified with
the Grassmannian ${\rm Gr}_{t}(t^{-n}k[[t]]^{d}/t^{n}k[[t]]^{d})$,
whose $K$-rational points, for $K$ a field containing $k$, are the
$K$-linear
subspaces of the finite dimensional $K$-vector space
$t^{-n}K[[t]]^{d}/t^{n}K[[t]]^{d}$ which are stable by multiplication by $t$;
in particular ${\cal Gr}_{(n)}$ is a projective variety. Furthermore, $\gr$
as a $k$-space is naturally isomorphic to the direct limit of the system of
projective varieties ${\cal Gr}_{(n)}$, hence $\gr$ is an ind-scheme. The
canonical morphism of $k$-spaces $\Theta :{\bf GL}_{d}(k((t)))\rightarrow
\gr
$ is a locally trivial fibration for the Zariski topology, {\it i.e.} $\gr$
is covered by open subsets over which $\Theta $ is a product.

\subsection{}\label{4.3} Denote by ${\rm M}_{d}$ the affine $k$-space of
$d$ by $d$ matrices. There
is canonical immersion of $k$-schemes $\iota :{\bf GL}_{d}(k[[t]])%
\hookrightarrow \cL ({\rm M}_{d}),$ which identifies ${\bf GL}_{d}(k[[t]])$
with the open subscheme defined by ${\rm det}\not= 0$. In particular we can
consider ${\bf GL}_{d}(k[[t]])$ as a semi-algebraic subset
of $\cL ({\rm M}_{d})$.
A subset of ${\bf GL}_{d}(k[[t]])$ will be called
semi-algebraic
(resp. $k [[t]]$-semi-algebraic)
if it is semi-algebraic
(resp. $k [[t]]$-semi-algebraic)
as a subset of $\cL ({\rm M}_{d})$. The space
${\bf GL}_{d}(k((t)))$ being an ind-scheme, one can associate to it an
underlying set, which by abuse we shall denote by the same letter. We shall
say a subset $A$ of ${\bf GL}_{d}(k((t)))$ is bounded $k [[t]]$-semi-algebraic
(resp. stable bounded $k [[t]]$-semi-algebraic, resp. cylindrical), if, for some
integer $n$ in ${\Bbb N}$, $t^{n}A$ is a $k [[t]]$-semi-algebraic
(resp.
stable $k [[t]]$-semi-algebraic, resp. cylindrical) subset of ${\bf
GL}_{d}(k[[t]])$.

We may extend the measures $\mu$ and $\tilde \mu$ to bounded $k [[t]]$%
-semi-algebraic, stable bounded $k [[t]]$-semi-algebraic and cylindrical subsets
respectively, by defining
\[
\mu (t^{-n} A) = \LL^{d^2n} \mu (A) \quad \text{and} \quad \tilde \mu
(t^{-n} A) = \LL^{d^2n} \tilde \mu (A),
\]
for $A$ a $k [[t]]$-semi-algebraic (resp. stable $k [[t]]$-semi-algebraic or
cylindrical) subset of ${\bf GL}_d(k [[t]])$, which is independent from the
choice of the integer $n$. One should also remark that the measures $\mu$
and $\tilde \mu$ are invariant under ${\bf GL}_d(k [[t]])$-action.

\subsection{}\label{4.4} For any integer $m$ in ${\Bbb Z}$, the functor
\[
R\longmapsto \Bigl\{M\in {\rm GL}_{d}(R((t))) : {\rm ord}_{t}{\rm det}
M=m\Bigr\}
\]
defines a subspace of ${\bf GL}_{d}(k((t)))$, which we shall denote by ${\bf %
GL}_{d}(k((t)))[m]$. Clearly
${\bf GL}_{d}(k((t)))[m]$ is an ind-scheme.

If $H$ is a linear subspace of $t^{-n}K[[t]]^{d}/t^{n}K[[t]]^{d}$, for
$K$ a field containing $k$,
we set
\begin{equation}
{\rm index}(H):=\dim t^{-n}K[[t]]^{d}/H-nd.  \label{index}
\end{equation}
We define ${\cal Gr}_{(n)}[m]$ as the projective variety parametrizing
subspaces in the Grassmannian ${\rm Gr}_{t}(t^{-n}k[[t]]^{d}/t^{n}k[[t]]^{d})
$ which are of index $m$. For fixed $m$, the varieties ${\cal Gr}_{(n)}[m]$
form an inductive system and we denote by $\gr [m]$ the corresponding
ind-scheme. For every $m$, the fibration $\Theta :{\bf GL}%
_{d}(k((t)))\rightarrow \gr$ restricts to a fibration $\Theta :{\bf GL}%
_{d}(k((t)))[m]\rightarrow \gr[m]$.

\subsection{}\label{4.5} Let $A$ be a subset of $\gr$. We say $A$ is
gr-stable at level $n$ if $A$ is
a constructible subset of ${\cal Gr}_{(n)}$. Note that in this case $A\cap %
\gr[m]$ is constructible in ${\cal Gr}_{(n)}[m]$ for every $m$. Hence we can
define the motivic measure $\tilde{\mu}_{\gr}(A)$ of $A$ as the element
\[
\tilde{\mu}_{\gr}(A):=\sum_{m\in {\Bbb Z}}\frac{[A\cap \gr[m]]}{\LL^{md}}
\]
in $\cM_{{\rm loc}}$, which makes sense, the sum on the right hand side
being finite.

Clearly if $A$ is gr-stable at level $n$ then $\Theta^{-1} (t^n A)$ is
cylindrical at level $2n$. The relation with the previously defined measure $%
\tilde \mu$ is given by the following Proposition.

\begin{proposition}
\label{comp}Let $A$ be a gr-stable subset of $\gr$. Then
\[
\tilde{\mu}_{\gr}(A)=(1-\LL^{-1})^{-1}\ldots (1-\LL^{-d})^{-1}\tilde{\mu}%
(\Theta ^{-1}(A)).
\]
\end{proposition}

\begin{proof}We may assume $A$ is contained in ${\cal Gr}_{(n)}[m]$. By Lemma
\ref{fib} applied to $r=m+nd$ and $p=2n$,
\[
\lbrack \pi _{2n}(\Theta ^{-1}(t^{n}A))]=[t^{n}A](\LL^{d}-1)\ldots (\LL^{d}-%
\LL^{d-1})\LL^{2nd^{2}-d(m+nd)},
\]
since $[{\rm GL}_{d, k}]=(\LL^{d}-1)\ldots (\LL^{d}-\LL^{d-1})$. We deduce
that
\[
\tilde{\mu}(\Theta ^{-1}(t^{n}A))=[t^{n}A](1-\LL^{-1})\ldots (1-\LL^{-d})\LL%
^{-nd^{2}}\LL^{-md}.
\]
The result follows since
\[
\tilde{\mu}(\Theta ^{-1}(t^{n}A))=\tilde{\mu}(t^{n}\Theta ^{-1}(A))=\LL%
^{-nd^{2}}\tilde{\mu}(\Theta ^{-1}(A))
\]
and $[t^{n}A]=[A]$.\end{proof}

\begin{lemma}
\label{fib}For any integer $p \geq 0$, the morphism
\[
\pi _{p}(\Theta ^{-1}({\rm Gr}_{t}(k[[t]]^{d}/t^{p}k[[t]]^{d})[r]))%
\longrightarrow {\rm Gr}_{t}(k[[t]]^{d}/t^{p}k[[t]]^{d})[r]
\]
is a locally trivial fibration for the Zariski topology with fiber
${\rm GL}_{d, k}
\times {\Bbb A}_{k}^{pd^{2}-dr}$. Here
${\rm GL}_{d, k}$ denotes
the algebraic variety of invertible
$d$ by $d$ matrices over
$k$.
\end{lemma}

\begin{proof}By taking
a cover of
the Grassmannian ${\rm Gr}_{t}(k[[t]]^{d}/t^{p}k[[t]]^{d})[r]$
by open Schubert cells corresponding to different choices of
bases of the lattice $k[[t]]^{d}$, one deduces the result from the
following elementary Lemma \ref{el}.\end{proof}

\begin{lemma}
\label{el}Let $p\geq m_{d}\geq \ldots \geq m_{1}\geq 0$ be integers.
Set $r=\sum_{1\leq i\leq d}m_{i}$.
Let $U$ be the subscheme of $\cL ({\rm M}_{d})$
consisting of triangular matrices $(a_{ij})$ with $a_{ii}=t^{m_{i}}$,
$a_{ij}=0$ for $j<i$, $m_i \leq {\rm ord}_t a_{ij}$ and
${\rm deg}_{t}a_{ij}<m_{j}$ for $j>i$. Consider
the morphism $\varphi~: U \times {\bf GL}_{d} (k[[t]])
\rightarrow U \times \cL_p ({\rm M}_{d})$
sending $(A, M)$  to $(A,  \pi_p (AM))$.
Then the image $W$ of $\varphi$
is an algebraic variety over $k$
which is isomorphic to $U\times {\rm GL}_{d, k}\times {\Bbb A}
_{k}^{pd^{2}-dr}$ by an isomorphism compatible with the projection $%
W\rightarrow U$.
\end{lemma}

Now we can express our
zeta functions
as integrals, defined as in \ref{int},
with
respect to the measure $\tilde{\mu}_{\gr}$
on the Grassmannian:
\begin{theorem}
\label{motivic=int on grass}Let ${\cal X}$ be a class of constructible $%
k[[t]]$-subalgebras. Then
\[
P_{L,{\cal X}}({\bf L}^{-s})=\int_{{\cal X\subset }{\rm Gr}_{0}^{d}(L)}{\bf L%
}^{{\rm index}(H)d-\left( \codim H\right) s}d\tilde{\mu}_{\gr}(H).
\]
\end{theorem}

\begin{proof} Let ${\cal X}_{n}$ be the subalgebras of codimension $n$. Then $%
{\cal X}_{n}\subset {\rm Gr}^{d}(t^{-n}L/t^{n}L)$ and
\begin{eqnarray*}
\left[ {\cal X}_{n}\right]  &=&\sum_{m\in {\Bbb Z}}\left[ {\cal X}_{n}\cap %
\gr[m]\right]  \\
&=&\sum_{m\in {\Bbb Z}}\frac{\left[ {\cal X}_{n}\cap \gr[m]\right] }{{\bf L}%
^{md}}{\bf L}^{md} \\
&=&\sum_{m\in {\Bbb Z}}\int_{{\cal X}_{n}\cap \gr[m]}{\bf L}^{{\rm index}%
(H)d}d\tilde{\mu}_{\gr}(H) \\
&=&\int_{{\cal X}_{n}}{\bf L}^{{\rm index}(H)d}d\tilde{\mu}_{\gr}(H).
\end{eqnarray*}
\end{proof}

\section{Rationality of motivic zeta functions of infinite dimensional Lie
algebras\label{Section Rationality}}

In this section we shall prove the following rationality result:

\begin{theorem}
\label{Rationality for k[[t]]}Let $k$ be a field of characteristic zero.

\begin{enumerate}
\item[(1)]  Let $L$ be a finite dimensional free $k[[t]]$-Lie algebra of the
form $L=L_{k}\otimes _{k}k[[t]]$, with $L_{k}$ a Lie algebra over $k$.
Let
${\cal X}_{t}^{\leq }$ (resp. ${\cal X}_{t}^{\triangleleft }$)
be the
class of subalgebras such that
${\cal X}_{t}^{\leq } (K)$ (resp. ${\cal X}_{t}^{\triangleleft } (K)$)
is the
set of all $K[[t]]$-subalgebras (resp. $K[[t]]$-ideals)
of $L \otimes K$, for every field
$K$ which is a finite extension of $k$.
Then $P_{L,{\cal X}_{t}^{*}}(T)$ for $*\in \{\leq ,\triangleleft \}$ is
rational,
belonging to ${\cal M}[T]_{{\rm loc}}$.

\item[(2)]  Let $L$ be a finite dimensional $k((t))$-Lie algebra and $L_{0}$
be a choice of some $k[[t]]$-Lie subalgebra of the form $L_{0}=L_{k}\otimes
_{k}k[[t]]$, with $L_{k}$ a Lie algebra over $k$.
Let
${\cal X}_{t,0}^{\leq }$ be the
class of subalgebras such that
${\cal X}_{t,0}^{\leq } (K)$
is the set of all $K [[t]]$-subalgebras $L \otimes K$
commensurable with $L_{0}
\otimes K$,
for every field
$K$ which is a finite extension of $k$.
Then $P_{L,{\cal X}_{0,t}^{\leq }}(T)$ is rational, belonging to ${\cal M}%
[T]_{{\rm loc}}$.
\end{enumerate}
\end{theorem}

Let $L$ be a finite dimensional free $k[[t]]$-Lie algebra of dimension $d$.
Choosing a basis allows us to identify $L$ additively with $k[[t]]^{d}$ and
we may view the ${\cal X}_{t}^{\leq }$ as a
subset of the Grassmannian $\gr$.

\begin{theorem}
\label{semi1}With the preceding notations, $\Theta ^{-1}({\cal X}_{t}^{*})$
is a $k [[t]]$-semi-algebraic subset of ${\bf GL}_{d}(k[[t]])$, for $*\in
\{\leq %
,\triangleleft \}$, and
\begin{eqnarray}
P_{L,{\cal X}_{t}^{*}}({\bf L}^{-s}) &=&\int_{{\cal X}_{t}^{*}}{\bf L}^{(d-s)%
\codim H}d\tilde{\mu}_{\gr}(H)  \label{666.1} \\
&=&(1-{\bf L}^{-1})^{-1}\ldots (1-{\bf L}^{-d})^{-1}\int_{\Theta ^{-1}({\cal X%
}_{t}^{*})}{\bf L}^{\left( d-s\right) {\rm ord}_{t}M}d\tilde{\mu}.
\label{666}
\end{eqnarray}
Furthermore, if $L=L_{k}\otimes _{k}k[[t]]$, with $L_{k}$ a Lie algebra over
$k$, then $\Theta ^{-1}({\cal X}_{t}^{*})$ is a semi-algebraic subset of $%
{\bf GL}_{d}(k[[t]])$.
\end{theorem}

\begin{proof} The first line (\ref{666.1}) follows from Theorem \ref
{motivic=int on grass} and the remark that $\codim H={\rm index}(H)$ if $H$
is a subalgebra of $L.$

Let $A^{*}=\Theta ^{-1}({\cal X}_{t}^{*}).$ Let $M=(m_{ij})$ be in ${\bf GL}%
_{d}(k[[t]])$ and write ${\bf m}_{i}=m_{i1}{e}_{1}+\cdots +m_{id}{e}_{d},$
for $i=1,\ldots ,d$. We have $M\in A^{\leq }$ if and only if, for every $1\leq
i,j\leq d$, there exist $Y_{ij}^{1},\ldots ,Y_{ij}^{d}\in K[[t]]$ such that $%
\beta ({\bf m}_{i},{\bf m}_{j})=\sum_{k=1}^{d}Y_{ij}^{k}{\bf m}_{k}$,
for some finite field extension $K$ of $k$. Here $%
\beta $ is the bilinear mapping $L\times L\rightarrow L$ corresponding to
the product in $L$. Let $C_{j}$ denote the matrix whose rows are ${\bf c}%
_{i}=\beta ({e}_{i},{e}_{j})$. We then have
\[
\beta ({\bf m}_{i},{\bf m}_{j})={\bf m}_{i}\left(
\sum_{l=j}^{d}m_{jl}C_{l}\right) .
\]
Then {\bf \ }$M\in A^{\leq }$ if and only if, for every $1\leq i,j\leq d$, one
can solve the matrix equation
\begin{equation}
{\bf m}_{i}\left( \sum_{l=j}^{d}m_{jl}C_{l}\right) =\left( Y_{ij}^{1},\ldots
,Y_{ij}^{d}\right) M  \label{matrix equation}
\end{equation}
for $Y_{ij}^{1},\ldots ,Y_{ij}^{d}\in K[[t]]$, with $K$ a finite field
extension of $k$. Let $M^{\prime }$ denote the
adjoint matrix, we can then rewrite (\ref{matrix equation}) as
\[
{\bf m}_{i}\left( \sum_{l=j}^{d}m_{jl}C_{l}\right) M^{\prime }=\left( \det
(M)Y_{ij}^{1},\ldots ,\det (M)Y_{ij}^{d}\right) .
\]
Let $g_{ijk}(m_{rs})$ denote the $k$-th entry of the $d$-tuple ${\bf m}
_{i}\left( \sum_{l=j}^{d}m_{jl}C_{l}\right) M^{\prime }$. Then the set $A^{%
\leq }$ has the following description:
\begin{equation}
A^{\leq }=\left\{ (m_{ij})\in \cL (M_{d}):{\rm ord}_{t}(\det (M))\leq {\rm %
ord}_{t}(g_{ijk}(m_{rs}))\text{ for }i,j,k\in \left\{ 1,\ldots ,d\right\}
\right\} .  \label{cone condition }
\end{equation}

The set $A^{\leq }$ is therefore $k [[t]]$-semi-algebraic. Let $\Phi ^{\leq
}(M)$
denote the conjunction of conditions ${\rm %
ord}_{t}(\det (M))\leq {\rm ord}_{t}(g_{ijk}(m_{rs})).$

Since ${\rm ord}_t: \cL (M_d)\rightarrow {\Bbb N\cup }\left\{ \infty
\right\} $ is a simple function this implies that $A_{n}$ is a
semi-algebraic set and in particular is constructible. As we promised in
section \ref{section: definition of motivic zeta} this provides an
alternative way to show the constructibility of this class of subalgebras.

We prove that $A^{\triangleleft }$ is $k [[t]]$-semi-algebraic. We
have
$M\in
A^{\triangleleft }$ if and only if, for every $1\leq i,j\leq d$, there exist $%
Y_{ij}^{1},\ldots ,Y_{ij}^{d}\in K[[t]]$, with $K$ some finite field
extension of $k$, such that $\beta ({\bf m}_{i},{\bf e%
}_{j})=\sum_{k=1}^{d}Y_{ij}^{k}{\bf m}_{k}.$ Let $g_{ijk}^{\triangleleft
}(m_{rs})$ denote the $k$-th entry of the $d$-tuple ${\bf m}%
_{i}C_{j}M^{\prime }$. Then an argument similar to the above implies that
the set $A^{\triangleleft }$ has the following description:
\[
A^{\triangleleft }=\left\{ (m_{ij})\in \cL (M_{d}):{\rm ord}_{t}(\det (M))%
\leq {\rm ord}_{t}(g_{ijk}^{\triangleleft }(m_{rs}))\text{ for }i,j,k\in
\left\{ 1,\ldots ,d\right\} \right\} .
\]
Hence $A^{\triangleleft }$ is $k [[t]]$-semi-algebraic. Indeed it is
defined by $%
\Phi ^{\triangleleft }$, the definable condition which is the conjunction of
conditions ${\rm ord}_{t}(\det (M))\leq {\rm ord}_{t}(g_{ijk}^{\triangleleft
}(m_{rs}))$. Equation (\ref{666}) follows from Proposition \ref{comp} and
the last statement is clear, since when $L$
is of the form $L_{k}\otimes _{k}k[[t]]$, the above definable
conditions are all semi-algebraic. \end{proof}

\begin{remark}
\label{conversion to triangular} The following provides an easier integral
expression in general to calculate. It is similar to
the integral expression that was
used in \cite{GSS}. Let $Y$ denote the variety ${\rm Tr}_{d, k}$ of $d\times d$
upper triangular matrices which is isomorphic to the affine space
${\Bbb A}_{k}^{d(d+1)/2}$.
Let $\tilde{\nu}$ denote the
measure on stable $k [[t]]$-semi-algebraic subsets of ${\cal L}(Y)$. For $%
i=1,\ldots ,d$ define the simple function $\alpha _{i}:{\cal L}%
(Y)\rightarrow {\Bbb Z}$ by $\alpha _{i}(M)={\rm ord}_{t}(m_{ii}).$ Then
\begin{eqnarray*}
&&(1-{\bf L}^{-1})^{-1}\ldots (1-{\bf L}^{-d})^{-1}\int_{A^{*}}{\bf L}%
^{\left( d-s\right) {\rm ord}_{t}(M)}d\tilde{\mu} \\
&=&(1-{\bf L}^{-1})^{-d}\int_{A^{*}\cap {\cal L}(Y)}{\bf L}^{-{\rm ord}%
_{t}(M)s}{\bf L}^{\alpha _{1}(M)+\cdots +i\alpha _{i}(M)+\cdots +d\alpha
_{d}(M)}d\tilde{\nu}.
\end{eqnarray*}
\end{remark}

\subsection{Commensurable subalgebras}\label{5.5} We turn now to the issue
of commensurable subalgebras. The following
discussion will also apply to counting commensurable subalgebras in a ${\Bbb %
Q}_{p}$-Lie algebra, an issue which has previously not been discussed.

We identify $L$ with $k((t))^{d}$ by a choice of a basis $e_{1},\ldots ,e_{d}
$ for the $k[[t]]$-submodule $L_{0}$ which is then identified with $%
k[[t]]^{d}$. The mapping $\Theta :{\bf GL}_{d}(k((t)))\rightarrow \gr$
sends a matrix $M$ to the lattice $\left\langle {\bf m}_{1},\ldots ,{\bf m}%
_{d}\right\rangle $ spanned by
\[
{\bf m}_{i}=m_{ii}e_{i}+\cdots +m_{id}e_{d},
\]
for $i=1,\ldots ,d$. The index then of $H=\left\langle {\bf m}_{1},\ldots ,%
{\bf m}_{d}\right\rangle $ (as defined in (\ref{index})) is just ${\rm ord}
_{t}(\det M).$

\begin{remark}
For subalgebras of $L_{0} \otimes K$ the index of $H$ is the same as the
codimension.
We might have considered encoding a zeta function of commensurable
subalgebras according to index but note that there are an infinite number of
commensurable submodules of index $0$ whilst there is only one of
codimension $0.$
\end{remark}

We now give a description of the codimension as a function of the entries of
$M.$ Recall that if the entries of $M$ are in the field $K$ containing
$k$,
the codimension of $H=\left\langle {\bf m}_{1},\ldots
,{\bf m}_{d}\right\rangle $ in $L \otimes K$ is defined as the sum of
the codimension
of $H \cap (L_{0} \otimes K)$ in $L_{0}\otimes K$ and in $H$:

\begin{lemma}
Let $M$ be a $K$-rational point of ${\bf GL}_{d}(k((t)))$ and define
\[
n_{i}=\min \left\{ 0,{\rm ord}_{t}(m_{ij}):j=1,\ldots ,d\right\} .
\]
Then the codimension of $H=\left\langle {\bf m}_{1},\ldots ,{\bf m}%
_{d}\right\rangle $ in $L \otimes K$ is the sum of ${\rm ord_{t}}(\det
(m_{ij}t^{n_{i}}))$ and $\sum_{i=1}^{d}n_{i}$.
\end{lemma}

\begin{proof}This follows from the fact that $H \cap (L_{0} \otimes K)
=\left\langle
t^{n_{1}}{\bf m}_{1},\ldots,t^{n_{d}}{\bf m}_{d}\right\rangle$. Hence the
codimension of $H \cap (L_{0} \otimes K)$ is given as usual by ${\rm
ord}_{t}(\det
(m_{ij}t^{n_{i}})).$ The codimension of $H \cap (L_{0} \otimes K)$ in $H$
on the other
hand is given by $\sum_{i=1}^{d}n_{i}$.\end{proof}

We define then the function
\[
\codim (M) :={\rm ord}_{t}(\det
(m_{ij}t^{n_{i}}))+\sum_{i=1}^{d}n_{i} =
{\rm ord}_{t}(\det
(m_{ij}))+ 2 \sum_{i=1}^{d}n_{i}
\] which is a simple function.

\subsection{} We want to convert the Poincar{{\'e}} series into integrals
over arc spaces.
To do this we think of $k((t))$ as consisting of two pieces: $k[[t]]$
identified with $k[[t]]$ via $v\mapsto v$ and $k((t))\backslash k[[t]]$
identified with $tk[[t]]$ via $v\mapsto v^{-1}$. Notice that the measure on $%
k((t))\backslash k[[t]]$ translates under this identification to the measure
${\bf L}^{2{\rm ord}_{t}(v)}\mu _{k[[t]]}$ on $tk[[t]].$ This follows
because the set $t^{-n}k[[t]]^{*}$ has measure ${\bf L}^{n}$ whilst $%
t^{n}k[[t]]^{*},$ its image under $v\mapsto v^{-1},$ has measure ${\bf L}%
^{-n}.$

We can partition ${\bf GL}_{d}(k((t)))$ into $2^{d^{2}}$ subsets ${\bf GL}
_{d}(k((t)))_{A}$ which we can index with an element of $A\in M_{d}(\left\{
\pm 1\right\} )$ so that
\[
{\bf GL}_{d}(k((t)))_{A}=\left\{ M\in {\bf GL}_{d}(k((t))):\ord_{t}(m_{ij})%
\geq 0\text{ if and only if }A_{ij}=1\right\} .
\]
We then identify each ${\bf GL}_{d}(k((t)))_{A}$ with a subset of $\cL ({\rm %
M}_{d})$ by the morphism $\Xi _{A}:M\mapsto \left(
m_{ij}^{A_{ij}}\right) $. Note that $\Xi _{A}\left( \Xi _{A}\left( M\right)
\right) =M.$

Now remark, that with the notations used in the proof of Proposition \ref
{semi1}, a matrix $M$ is in $\Theta ^{-1}({\cal X}_{0}^{\leq })$ if and only
if $\Phi ^{\leq }(M)$ is true where $\Phi ^{\leq }$ was defined as the
conjunction of statements ${\rm ord}_{t}(\det (M))\leq {\rm ord}%
_{t}(g_{ijk}(m_{rs}))$ appearing in (\ref{cone condition }). Hence the
subset $\Xi _{A}\left( \Theta ^{-1}({\cal X}_{t,0}^{\leq })\cap {\bf GL}%
_{d}(k((t)))_{A}\right) $ is a $k [[t]]$-semi-algebraic subset
of
$\cL ({\rm M}_{d})$ defined by the condition
$\Phi ^{\leq }(\Xi _{A}M)$ is true.
Furthermore, when $L_0$
is of the form $L_{k}\otimes_{k}k[[t]]$, the definable
conditions considered are all semi-algebraic, hence
$\Xi _{A}\left( \Theta ^{-1}({\cal X}_{t,0}^{\leq })\cap {\bf GL}%
_{d}(k((t)))_{A}\right) $ is  semi-algebraic in this case.

Hence we have the following result:

\begin{theorem}
\label{semi2}With the preceding notations, $\Xi _{A}\left( \Theta ^{-1}(%
{\cal X}_{t,0}^{\leq })\cap {\bf GL}_{d}(k((t)))_{A}\right) $ is a $k [[t]]$%
-semi-algebraic subset of $\cL ({\rm M}_{d})$, and

\begin{eqnarray}
P_{L,{\cal X}_{t,0}^{\leq }}({\bf L}^{-s})&=&\int_{\Theta ^{-1}({\cal X}%
_{t,0}^{\leq })}{\bf L}^{{\rm ord}_{t}(\det M)d-\left( \codim M\right) s}d%
\tilde{\mu}  \nonumber \\
&=&\sum_{A\in {M}_{d}(\left\{ \pm 1\right\} )}\int_{\Xi _{A}\left( \Theta
^{-1}({\cal X}_{t,0}^{\leq })\cap {\bf GL}_{d}(k((t)))_{A}\right) }{\bf L}^{%
{\rm ord}_{t}(\det \Xi _{A}\left( M\right) )d-\left( \codim \Xi _{A}\left(
M\right) \right) s}  \label{777} \\
&&\times {\bf L}^{\sum_{A_{ij}=-1}2{\rm ord}_{t}m_{ij}^{-A_{ij}}}d\tilde{\mu}.
\nonumber
\end{eqnarray}

Furthermore, if $L_{0}=L_{k}\otimes_{k}k[[t]]$, with $L_{k}$ a Lie algebra
over $k$, then
\[
\Xi _{A}\left( \Theta ^{-1}({\cal X}_{t,0}^{\leq })\cap {\bf GL}%
_{d}(k((t)))_{A}\right)
\]
is a semi-algebraic subset of $\cL (M_{d})$.
\end{theorem}

\begin{proof}Everything has already been proven, except for (\ref{777}),
which follows from Proposition \ref{comp}. \end{proof}

\subsection{Proof of Theorem \ref{Rationality for k[[t]]}}

Let us first prove (1). Since, by Theorem \ref{semi1}, we have (\ref{666})
the result follows from Theorem \ref{5.1prime}, Remark \ref{min} and the
last statement in Theorem \ref{semi1}. Similarly to deduce (2) from Theorem
\ref{semi2}, Theorem \ref{5.1prime} and Remark \ref{min}, one needs to check
that there exists an integer $N$ such that
\[
N\codim \Xi _{A}\left( M\right) >{\rm ord}_{t}(\det \Xi _{A}\left( M\right)
)d+\sum_{A_{ij}=-1}2{\rm ord}_{t}m_{ij}^{-A_{ij}}.
\]
This follows from the fact that $\codim \Xi _{A}\left( M\right) >{\rm ord}%
_{t}(\det \Xi _{A}\left( M\right) )$ and if $A_{ij}=-1$ then ${\rm ord}%
_{t}m_{ij}^{-1}<n_{i}.$ Hence if we take $N=3d$ we are done. \hfill \qed

\begin{remark}
If Theorem \ref{5.1prime} remains true when semi-algebraic is replaced
everywhere by $k [[t]]$-semi-algebraic and simple by $k [[t]]$-simple, as
suggested in
Remark \ref{bad}, then Theorem \ref{Rationality for k[[t]]} remains true
with the same proof, without assuming that $L$ is obtained by extension of
scalars from a Lie algebra $L_{k}$ defined over $k$.
\end{remark}

\subsection{Commensurable subalgebras in $p$-adic Lie algebras}
Note that the analysis above can also be applied to any finite dimensional $%
{\Bbb Q}_{p}$-Lie algebra ${\cal L}$ with a choice of ${\Bbb Z}_{p}$-lattice
$L$ where we define the commensurable zeta function as:
\[
\zeta _{{\cal L},L}(s)=\sum_{H\in {\cal X}({\cal L},L)}|L+H:L\cap H|^{-s}
\]
where ${\cal X}({\cal L},L)$ denotes theset of  ${\Bbb Z}_{p}$-subalgebras
commensurable with $L$. This is the analogue of what we did above where the
index $|L+H:L\cap H|$ corresponds to the codimension since $\codim H=\codim
(L_{0}+H:L\cap H)$ in the $k[[t]]$-setting. We shall consider later, see
Corollary \ref{np}, the question of zeta functions of $L\otimes {\Bbb Z}_{p},
$ where $L$ is ${\Bbb Z}$-Lie algebra, but here the analysis applies to a
Lie algebra with no underlying ${\Bbb Z}$-structure. So it is worth
recording the following result for $\zeta _{{\cal L},L}(s),$ proved as we
said by the same analysis as above replacing ${\bf L}$ by $p$, $k[[t]]$ by $%
{\Bbb Z}_{p}$, $k((t))$ by ${\Bbb Q}_{p}$, and Theorem \ref{5.1prime} by
Theorem 3.2 of \cite{rat}.

\begin{theorem}
Let ${\cal L}$ be a finite dimensional ${\Bbb Q}_{p}$-Lie algebra with a
choice of ${\Bbb Z}_{p}$-lattice $L$ such that $L\otimes {\Bbb Q}_{p}={\cal L%
}.$ Then $\zeta _{{\cal L},L}(s)$ is a rational function in $p^{-s}.$
\end{theorem}

\subsection{${\Bbb F}_{p}[[t]]$-Lie algebras} We mention that there is
another open rationality question which fits into
this context. This concerns Lie algebras in characteristic $p.$ For example
the following is conjectured in section 5 of \cite{dus-sl2}:

\begin{conjecture}
Let $L$ be a finite dimensional Lie algebra over ${\Bbb F}_{p}[[t]]$
additively isomorphic to ${\Bbb F}_{p}[[t]]^{d}.$ Then
\[
\zeta _{L}(s)=\sum_{H\leq L}|L:H|^{-s}
\]
is a rational function in $p^{-s}$ where the sum is taken over ${\Bbb F}
_{p}[[t]]$-subalgebras $H.$
\end{conjecture}

\subsection{}The techniques involved in proving such a result in
characteristic zero (quantifier elimination and
Hironaka's resolution of singularities) do not exist at present in
characteristic $p$. However, in \cite{dus-sl2} the rationality was proved
for one class of ${\Bbb F}_{p}[[t]]$-algebras. Take $L$ to be an algebra
additively isomorphic to ${\Bbb Z}^{d}$ and set $L_{p}=L\otimes {\Bbb F}%
_{p}[[t]],$ then one can prove that for almost all primes $\zeta _{L_{p}}(s)$
is a rational function in $p^{-s}.$ In fact one gets that $\zeta
_{L_{p}}(s)=\zeta _{L\otimes {\Bbb Z}_{p}}(s)$ for almost all primes. The
proof actually follows the theme of this current paper, that there is an
integral over the additive Haar measure on $\left( {\Bbb F}_{p}[[t]]\right)
^{d}$ representing this zeta function ${\Bbb F}_{p}[[t]]$ which formally
looks the same as the integral representing $\zeta _{L\otimes {\Bbb Z}%
_{p}}(s).$ Macintyre then observed that the calculation of this latter
integral will carry over formally to the former essentially for any primes
for which one does not divide by in this calculation. See also Theorem 8.3.2
of \cite{def} for much more general results along these lines.

\subsection{}If we are going to consider the possible rationality
of some of the other
motivic zeta functions considered in this paper then it is necessary to
introduce the concept of {\em polynomial growth} for $(L,{\cal X})$. This is
going to be a necessary requirement if we are to prove that the zeta
function $P_{L,{\cal X}}(T)$ can be expressed as an element of ${\cal M}[T]_{%
{\rm loc}}$.

\begin{definition}
We say that $(L,{\cal X})$ has {\em polynomial growth} if there exists some $%
d\in {\Bbb N}$ such, that for all $n$, $A_{n}({\cal X})$ is a constructible
set and
\[
\dim A_{n}({\cal X})\leq dn.
\]
Recall
a constructible set is of dimension $\leq m$, if it may be partioned into
finitely many locally closed pieces with dimension $\leq m$.
\end{definition}

It is a necessary condition that $(L,{\cal X})$ have polynomial growth for
the zeta function to be rational, because the coefficients of a rational
series satisfy a linear recurrences of finite order. For example, when $L$
is a finite dimensional $k[[t]]$-algebra additively isomorphic to $%
k[[t]]^{d} $ and ${\cal X}$ the class of commensurable $k[[t]]$-subalgebras,
then $A_{n}({\cal X})$
is contained in ${\cal Gr}_{(n)}$ whose dimension is
bounded by
$2nd^{2}$, the dimension of ${\rm M}_{d}(t^{-n}k[[t]]/t^{n}k[[t]])$.

\begin{problem}Characterize the ${\Bbb N}$-filtered or ${\Bbb Z}$-filtered
infinite
dimensional Lie algebras for which ${\cal X}^{\leq }$ has polynomial growth.
\end{problem}

\subsection{Polynomial growth} If we consider the algebra ${\frak
d}^{+}:=\bigoplus_{j\geq 1}{k}\cdot d_{j}$,
where $(d_{i},d_{j})=\left( i-j\right) d_{i+j}$ and ${\cal X}$ is
the class of
all subalgebras, then this does not have polynomial growth. However, if we
define ${\cal X}^{d}$ to be the class of
$d$-generated subalgebras, we find that $(%
{\frak d}^{+},{\cal X}^{d})$ does have polynomial growth.

\begin{problem}
Calculate $P_{{\frak d}^{+},{\cal X}^{d}}(T)$ for $d\geq 2.$ Is it rational?
\end{problem}

More generally, $(L,{\cal X}^{d})$ has polynomial growth if $L$ is a
well-covered Lie algebra (where
recall a subalgebra of codimension $n$ in $L \otimes K$
contains $L_{f(n)} \otimes K$ for some function $f(n)$)
and $\dim (L/L_{f(n)})\leq
n^{c}$ for some fixed $c$. This is true for example for all the simple
graded Lie algebras of finite growth. However we are unable to prove that
the corresponding zeta function encoding these $d$-generated subalgebras is
rational.

\section{$k [[t]]$-powered nilpotent groups}

\subsection{}In \cite{GSS}
subgroups of finite index in a torsion-free finitely generated
nilpotent group $G$ were encoded in zeta functions. It was shown there
how,
for almost all primes $p$, there is a one-to-one correspondence between
subgroups of finite index in the pro-$p$ completion $\widehat{G}_{p}$ of $G$
and subalgebras of finite index in $L\otimes {\Bbb Z}_{p}$ for an associated
Lie algebra $L$ over ${\Bbb Z}.$ The pro-$p$ completion has the structure of
a ${\Bbb Z}_{p}$-powered nilpotent group.

We can use the same ideas to show that the motivic zeta functions counting
subalgebras in a finite dimensional $k[[t]]$-Lie algebra are also encoding
the subgroup structure of a $k[[t]]$-powered nilpotent group $G$ where
subgroups are taken to be $k[[t]]$-closed subgroups.

An $R$-powered nilpotent group is a nilpotent group with an action of the
ring $R$ satisfying various axioms which can be found in Chapter 6 of P.
Hall's notes on nilpotent groups \cite{Hall}.

We shall be interested in a certain type of $R$-powered nilpotent group
which arises when $R$ is a binomial ring in the following way. Take a
torsion-free finitely generated nilpotent group $G$ with a choice of Malcev
basis $x_{1},\ldots ,x_{d}.$ The group $G$ can then be identified with the
set of elements ${\bf x}({\bf a})=x_{1}^{a_{1}}\ldots x_{d}^{a_{d}}$ with $%
{\bf a}=(a_{1},\ldots ,a_{d})\in {\Bbb Z}^{d}.$ In \cite{Hall} P. Hall
proved that there exist polynomials over ${\Bbb Q}$ in suitably many
variables which define mappings
\begin{eqnarray*}
\lambda &:&{\Bbb Z}^{d}\times {\Bbb Z\rightarrow Z} \\
\mu &:&{\Bbb Z}^{d}\times {\Bbb Z}^{d}\rightarrow {\Bbb Z}^{d}
\end{eqnarray*}
such that for ${\bf a}_{1}$ and ${\bf a}_{2}\in {\Bbb Z}^{d}$ and $k\in
{\Bbb Z}$%
\begin{eqnarray*}
{\bf x}({\bf a})^{k} &=&{\bf x(}\lambda ({\bf a},k)) \\
{\bf x}({\bf a}_{1})\cdot {\bf x}({\bf a}_{2}) &=&{\bf x}(\mu ({\bf a}_{1},%
{\bf a}_{2})).
\end{eqnarray*}
Since these polynomials over ${\Bbb Q}$ take integer values, they can be
represented as ${\Bbb Z}$-linear combinations of binomial polynomials.
Suppose $R$ is a binomial ring, {\it i.e.} $\binom{r}{n}$ is a well-defined
element of $R$ for $r\in R$ and $n\in {\Bbb Z}$. Then we can use the
polynomials $\lambda $ and $\mu $ to define the structure of an $R$-powered
nilpotent group on the set of elements ${\bf x}({\bf a})=x_{1}^{a_{1}}\ldots
x_{d}^{a_{d}}$ with ${\bf a}=(a_{1},\ldots ,a_{d})\in R^{d}.$

Therefore to every finitely generated torsion-free nilpotent group $G$ we
can define a $k[[t]]$-powered group which we shall denote $G^{k[[t]]}.$
There exists a ${\Bbb Q}$-rational representation of $G^{k[[t]]}$ embedding
it as a subgroup of ${\rm Tr}_{n}^{1}(k[[t]])$, the upper triangular matrix
group of dimension $n$ with $1$'s on the diagonal. The image of $G^{k[[t]]}$
under the map $\log $ defines a nilpotent $k[[t]]$-subalgebra ${\cal L}=\log
G^{k[[t]]}$ of ${\rm Tr}_{n}^{0}(k[[t]])$, the upper triangular matrix group
of dimension $n$ with $0$'s on the diagonal. (Note that $\log $ is a
polynomial map since $\left( U-1\right) ^{n}=0$ if $U\in {\rm Tr}
_{n}^{1}(k[[t]]).)$ The Lie algebra ${\cal L}$ is $L_{{\Bbb Q}}\otimes
k[[t]] $ where $L_{{\Bbb Q}}$ is the Lie algebra corresponding to $G^{{\Bbb Q%
}}$ under the Malcev correspondence. In general the image of $G$ under $\log
$ is not a Lie algebra. However there is some integer $f$ such that $L=\log
G^{f}$ is a ${\Bbb Z}$-Lie algebra. This is the Lie algebra mentioned in the
first paragraph used by Grunewald, Segal and Smith to set up a one-to-one
correspondence between subgroups in $G^{{\Bbb Z}_{p}}=\widehat{G}_{p}$ and $%
L\otimes {\Bbb Z}_{p}$ for all primes not dividing $f$ (see Theorem 4.1 of
\cite{GSS}).

A subgroup $H$ of $G^{k[[t]]}$ is $k[[t]]$-closed if $h^{r}\in H$ for all $%
r\in k[[t]]$ and $h\in H.$ We shall define the codimension of such a
subgroup as the length $n$ of the maximal chain of $k[[t]]$-closed subgroups
$H_{i}$ $0\leq i\leq n$ with $H_{i}$ contained in $H_{i-1}$ with infinite
index and $H_{0}=G^{k[[t]]}$ and $H_{n}=H.$

The claim is then that $\log $ gives a way of getting from $k[[t]]$-closed
subgroups to $k[[t]]$-subalgebras to prove the following:

\begin{theorem}
\label{one-to-one corr}There is a one-to-one correspondence between the set
of $k[[t]]$-closed subgroups of $G^{k[[t]]}$ of finite codimension and the
set of $k[[t]]$-subalgebras of finite codimension in ${\cal L}$ which
preserves the codimension.
\end{theorem}

\begin{proof} We need to prove that the image of a $k[[t]]$-closed subgroup $%
H $ under the map $\log $ is a $k[[t]]$-subalgebra. The proof of Lemma 1 in
Chapter 6 of \cite{Segal-polycyclic groups} can be applied in this setting
to prove that $m{\Bbb Z}\log H\subset \log H.$ But since $H$ is
$k[[t]]$-closed and $k[[t]]\supset {\Bbb Q}$ we have $h^{1/m}\in H.$ Hence
${\Bbb Z}%
\log H\subset \log H,$ {\it i.e. }$\log H$ is additively closed. Once we
have this then Corollary 3 of Chapter 6 of \cite{Segal-polycyclic groups}
implies that $\log H$ is also closed under the Lie bracket by applying the
inverse of the Campbell-Hausdorff series. We have implicitly used the
identity $r\log h=\log h^{r}$ above when $r\in {\Bbb Q}$. This is true in
fact for all values of $r\in k[[t]].$ This follows from the fact that there
exist polynomials $f_{1}({\bf a}),\ldots ,f_{d}({\bf a})$ such that $\log
x_{1}^{a_{1}}\ldots x_{d}^{a_{d}}=f_{1}({\bf a})\log x_{1}+\cdots +f_{d}(%
{\bf a})\log x_{d}.$ Since $r\log h=\log h^{r}$ is true for all $h\in G$ and
$r\in {\Bbb Z}$ we get that the following is a polynomial identity for all $%
{\bf a}\in k[[t]]^{d}$ and $r\in {\Bbb Z}$%
\[
f_{i}(\lambda ({\bf a},r))=rf_{i}({\bf a}).
\]
Hence this must be a formal identity of polynomials which is then true for
all values of $r\in k[[t]].$ Conversely, if $M$ is a $k[[t]]$-subalgebra of $%
{\cal L}$ then $\exp M$ is well-defined on ${\rm Tr}_{n}^{0}(k[[t]])$ being
a finite polynomial map. By using the
Campbell-Hausdorff series one obtains directly that
$\exp M$ is a subgroup and the identity $r\log h=\log h^{r}$ again implies
that $\exp M$ is a $k[[t]]$-closed subgroup taking $h=\exp m$ and using the
identity $\log \exp (m)=m.$ Since $\log $ and $\exp $ are bijective maps, we
have a one-to-one correspondence as detailed in the statement of the Theorem.
The correspondence preserves codimension since this is defined by lengths of
chains of $k[[t]]$-closed subgroups or $k[[t]]$-subalgebras. \end{proof}

We can use this to deduce the existence and rationality of a motivic zeta
function associated to $G^{k[[t]]},$ a certain type of loop group.

\begin{corollary}
Let $G$ be a finitely generated torsion-free nilpotent group and $G^{k[[t]]}$
the associated $k[[t]]$-powered nilpotent group.
For any field $K$ which is a finite extension of $k$, let
${\cal X} (K)$ denote the
set of $K[[t]]$-closed subgroups of $G^{K[[t]]}$ and
let $A_n (K)$
denote the
set of $K[[t]]$-closed subgroups of $G^{K[[t]]}$ of codimension $n$.
There is a natural structure of constructible set on $A_n = \cup_K A_n
(K)$,
hence
\[
P_{G^{k[[t]]},{\cal X}}(T) :=\sum_{n \in {\Bbb N}} \, [A_n] \, T^n
\]
is a well-defined element of ${\cal M}[[T]].$ We call it the
{\em motivic
zeta function of }$G^{k[[t]]}.$ It is equal to the motivic zeta function $P_{%
{\cal L},{\cal X}_{t}^{\leq}}(T)$ of the associated $k[[t]]$-Lie algebra
${\cal L}%
=\log G^{k[[t]]}=\log G^{{\Bbb Q}}\otimes _{{\Bbb Q}}k[[t]]$. Hence
$P_{G^{k[[t]]},{\cal X}}(T)$ is rational, belonging to ${\cal M}[T]_{{\rm
loc}}$.
\end{corollary}

We have a similar result for normal $k[[t]]$-closed subgroups in $G^{k[[t]]}$
since the one-to-one correspondence in Theorem \ref{one-to-one corr} sends
normal subgroups to ideals.

\section{Motivic and $p$-adic cone integrals\label%
{Section:Motivic cone integrals}}

The motivic integrals we have defined to capture the subalgebras in a $%
k[[t]] $-algebra are examples of integrals that the first author and
Grunewald studied in the $p$-adic setting in \cite{duSG-Abscissa} called
cone integrals. We make the same definition in the motivic setting:

\begin{definition}
We call an integral
\[
Z_{{\cal D},{\rm geom}}(s)
:= \int_{V}{\bf L}^{-{\rm ord}_{t}(f_{0})s-{\rm ord}
_{t}(g_{0})}d\mu
\]
{\em a motivic cone integral} if $f_{0}({\bf x})$ and $g_{0}({\bf x})$ are
polynomials in $k[[t]][x_{1},\ldots ,x_{m}]$ and there exist polynomials $%
f_{i}({\bf x}),g_{i}({\bf x})$, $i=1,\ldots ,l$, in $k[[t]][x_{1},\ldots
,x_{m}]$ such that
\[
V=\left\{ {\bf x}\in \cL ({\Bbb A}_{k}^{m}):{\rm ord}_{t}(f_{i}({\bf x)})%
\leq {\rm ord}_{t}(g_{i}({\bf x}))\text{ for }i=1,\ldots ,l\right\} .
\]
The set ${\cal D}=\left\{ f_{0},g_{0},f_{1},g_{1},\ldots
,f_{l},g_{l}\right\} $ is called the {\em set of cone integral data}{\it .}
\end{definition}

\begin{definition}
Let $L$ be a $k[[t]]$-algebra additively isomorphic to $\left( k[[t]]\right)
^{d}$ (or a ${\Bbb Z}$-Lie algebra additively isomorphic to ${\Bbb Z}^{d}).$
Then we define $Z_{L,{\rm geom}}(s)$, {\em the motivic cone integral
associated to }$L$, to be
\begin{eqnarray*}
Z_{L,{\rm geom}}(s) &:=&\int_{A}{\bf L}^{-s{\rm ord}_{t}(M)}d\mu  \\
&=&Z_{{\cal D},{\rm geom}}(s)
\end{eqnarray*}
where
\[
A=\left\{ (m_{ij})\in \cL (M_{d}):{\rm ord}_{t}(\det (M))\leq {\rm ord}
_{t}(g_{ijk}(m_{rs}))\text{ for }i,j,k\in \left\{ 1,\ldots ,d\right\}
\right\}
\]
and ${\cal D}$ is the corresponding cone data.
\end{definition}

In section \ref{Section Rationality} we therefore established that
\[
P_{L,{\cal X}_{t}^{\leq }}({\bf L}^{-s})=(1-{\bf L}^{-1})^{-1}
\ldots (1-{\bf L}
^{-d})^{-1} Z_{L,{\rm geom}}(s-d).
\]

The set $V$ in the definition of the motivic cone integral is $k [[t]]$%
-semi-algebraic and therefore measurable; the functions ${\rm ord}_{t}(f_{0})
$ and ${\rm ord}_{t}(g_{0})$ are $k [[t]]$-simple functions. However, in
this more
general setting, we don't know whether the fibres of ${\rm ord}
_{t}(f_{0}):V\rightarrow {\Bbb N}$ and ${\rm ord}_{t}(g_{0}):V\rightarrow
{\Bbb N}$ are weakly stable. This is why we need to use the completion $%
\widehat{{\cal M}}$ of ${\cal M}_{{\rm loc}}$ with respect to the filtration
$F^{m}{\cal M}_{{\rm loc}}$, introduced in \ref{completion}, and to replace
the measure $\tilde{\mu}$ by the measure $\mu $. Hence, Theorem \ref{5.1}
gives us the following:

\begin{theorem}
\label{rationality of motivic cone}Assume the polynomials $f_{i}({\bf x}%
),g_{i}({\bf x})$, $i=0,\ldots ,l,$ have their coefficients in $k$. The
motivic cone integral $Z_{{\cal D},{\rm geom}}(s)$ belongs to the subring $%
\overline{\cN}$ of the ring $\widehat{{\cal M}}[[T]]$ which is generated by
the image in $\widehat{{\cal M}}[[T]]$ of ${\cal M}_{{\rm loc}}[T]$, $({\bf L%
}^{i}-1)^{-1}$ and $(1-{\bf L}^{-a}T^{b})^{-1}$ with $i\in {\Bbb N}$, $a\in
{\Bbb N}$, $b\in {\Bbb N}\setminus \{0\}$, where $T={\bf L}^{-s}$.
\end{theorem}

In the algebra setting we had that the cone conditions were in fact stable
and hence, by Theorem \ref{5.1prime}, the associated zeta function $Z_{L,%
{\rm geom}}(s)$ is an element of the subring of ${\cal M}_{{\rm loc}}[[T]]$,
where $T={\bf L}^{-s}.$ However in the more general setting here, we are
required to take infinite sums of elements of ${\cal M}_{{\rm loc}}$ which
are then only defined in the associated complete ring $\widehat{{\cal M}}.$
At the moment it is not known if the image of ${\cal M}_{{\rm loc}}$ is the
same as ${\cal M}_{{\rm loc}}$ since there might be some kernel. The
infinite sums of elements of ${\cal M}_{{\rm loc}}$ and the terms like $(%
{\bf L}^{i}-1)^{-1}$ which don't occur in Theorem \ref{5.1prime} appear for
example when we consider the constant term of these motivic zeta functions
which take the following form:
\[
\int_{V\cap \{{\bf x} : \ord_t (f_{0}({\bf x})) = 0\}} {\bf L}^{-{\rm ord}%
_{t}(g_{0})}d\mu .
\]
Even if $g_{0}=0$ we get terms coming from the cone conditions like $({\bf L}%
^{i}-1)^{-1}$ as we shall see when we derive an explicit expression for
motivic cone integrals.

In \cite{duSG-Abscissa} the first author and Grunewald produced an explicit
formula for $p$-adic cone integrals extending Denef's explicit formula for
the Igusa zeta function. We can do the same thing for the motivic cone
integrals. As we shall explain, this will impact on how canonical the
explicit formula for $p$-adic cone integrals is. It will also allow us to
define a topological zeta function associated to $p$-adic cone integrals and
ultimately to the zeta function of an algebra defined over ${\Bbb Z}.$

\subsection{}Assume the polynomials $f_{i}({\bf x}),g_{i}({\bf x})$,
$i=0,\ldots,l,$ have
their coefficients in $k$. Let $F=f_{0}g_{0} \ldots f_{l}g_{l}$ and denote
by $D$ the divisor defined by $F = 0$. Let $(Y,h)$ be a resolution of $F$.
By this, we mean that $Y$ is a smooth and connected algebraic variety, $h :
Y \rightarrow X = {\Bbb A}_k^{m}$ is proper, that the restriction $h : Y
\setminus h^{-1} (D) \rightarrow X \setminus D$ is an isomorphism, and that $
(h^{-1} (D))_{{\rm red}}$ has only normal crossings as a subvariety of $Y$.
Let $E_{i},$ $i\in T$ be the irreducible (smooth) components of $%
(h^{-1}(D))_{{\rm red}}$ where $D$ is the divisor defined by $F=0$ in $X=%
{\Bbb A}^{m}.$ For each $i\in T,$ denote by $N_{i}(f_{j})$ and $N_{i}(g_{j})$
the multiplicity of $E_{i}$ in the divisor of $f_{j}\circ h$ and $g_{j}\circ
h$ on $Y$ and by $\nu_{i}-1$ the multiplicity of $E_{i}$ in the divisor of $%
h^{*}dx$, where $dx$ is a local non vanishing volume form. For $i\in T$ and $%
I\subset T,$ we consider the schemes $E_{i}^{\circ }=E_{i}\backslash
\bigcup_{j\neq i}E_{j},E_{I}=\bigcap_{i\in I}E_{i}$ and $E_{I}^{\circ
}=E_{I}\backslash \bigcup_{j\in T\backslash I}E_{j}.$ When $I=\emptyset ,$
we have $E_{\emptyset }=Y.$

Define a closed cone $\overline{D_{T}}$ as follows:
\[
\overline{D_{T}}
:=\left\{ (x_{1},\ldots,x_{t})\in {\Bbb R}_{\geq
0}^{t}:\sum_{j=1}^{t}N_{j}(f_{i})x_{j}\leq \sum_{j=1}^{t}N_{j}(g_{i})x_{j}%
\text{ for }i=1,\ldots,l\right\},
\]
where ${\rm card}T=t$ and ${\Bbb R}_{\geq 0}=\left\{ x\in {\Bbb R}:x\geq
0\right\} $. Denote the lattice points in $\overline{D_{T}}$ by $\overline{%
\Delta_{T}},$ {\it i.e.} $\overline{\Delta_{T}}=\overline{D_{T}}\cap {\Bbb N}%
^{t}.$ We can write this cone as a disjoint union of open simplicial pieces
that we shall call $R_{k}$, $k=0,1,\ldots,w.$ We shall assume that $%
R_{0}=(0,\ldots,0)$ and that the next $q$ pieces are all the open
one-dimensional edges in the cone $\overline{D_{T}}$: for $k=1,\ldots,q$,
\[
R_{k}=\left\{ \alpha {\bf e}_{k}=\alpha (q_{k1},\ldots,q_{kt}):\alpha
>0\right\} .
\]

Since these are all the edges, for any $k\in \left\{ 0,\ldots,w\right\} $
there exists some subset $M_{k}\subset \left\{ 1,\ldots,q\right\} $ such
that
\[
R_{k}=\left\{ \sum_{j\in M_{k}}\alpha_{j}{\bf e}_{j}:\alpha_{j}>0\text{ for
all }j\in M_{k}\right\} .
\]
Note that $m_{k}:={\rm card}M_{k}\leq t.$ We can choose the decomposition of
the cone $\overline{D_{T}}$ so that for any $k\in \left\{ 0,\ldots,w\right\}
$%
\[
R_{k}\cap {\Bbb N}^{t}=\left\{ \sum_{j\in M_{k}}l_{j}{\bf e}_{j}:l_{j}\in
{\Bbb N}_{>0}\text{ for all }j\in M_{k}\right\},
\]
as follows from \cite{Danilov}, p.123-124.

Define for each $k=1,\ldots,q$ the following non-negative constants:
\begin{eqnarray}
A_{k} &=&\sum_{j=1}^{t}q_{kj}N_{j}(f_{0})  \label{A_kB_k} \\
B_{k} &=&\sum_{j=1}^{t}q_{kj}\left( N_{j}(g_{0})+\nu_{j}\right) .  \nonumber
\end{eqnarray}

For each $I\subset T$ define:
\begin{eqnarray*}
D_{I} &=&\left\{ \left( k_{1},\ldots,k_{t}\right) \in \overline{D_{T}}%
:k_{i}>0\text{ if }i\in I\text{ and }k_{i}=0\text{ if }i\in T\backslash
I\right\} \\
\Delta_{I} &=&D_{I}\cap {\Bbb N}^{t}.
\end{eqnarray*}

So $\overline{D_{T}}=\bigcup_{I\subset T}D_{I}$, a disjoint union.

For each $I\subset T$ there is then a subset $W_{I}\subset \left\{
0,\ldots,w\right\} $ so that
\[
D_{I}=\bigcup_{k\in W_{I}}R_{k}.
\]
If $k \in W_I$, we set $I_k = I$.
We can now state the explicit formula for motivic cone integrals:

\begin{theorem}
\label{explicit expression for motivic cone}Assume the polynomials $f_{i}(%
{\bf x}),g_{i}({\bf x})$, $i=0,\ldots ,l,$ have their coefficients in $k$.
The following equality holds in $\overline{\cN}$:
\begin{equation}
Z_{{\cal D},{\rm geom}}(s)=\sum_{k=0}^{w}({\bf L}-1)^{|I_{k}|}{\bf L}%
^{-m}[E_{I_{k}}^{\circ }]\prod_{j\in M_{k}}\frac{{\bf L}^{-(A_{j}s+B_{j})}}{%
(1-{\bf L}^{-(A_{j}s+B_{j})})}.  \label{form}
\end{equation}
\end{theorem}

\begin{proof}The proof works just the same as the one in \cite{duSG-Abscissa}
for the explicit formula for $p$-adic cone integrals, using Lemma 3.4 of
\cite{Denef-Loeser-Arcs} for performing the change of variable, cf. the
proof of Proposition 2.2.2 of \cite{Denef-Loeser-motivic Igusa}. \end{proof}

\subsection{}\label{over}Our goal now is to explain the relationship
between motivic
cone integrals
and $p$-adic cone integrals. We begin with some preliminaries. We assume $k=%
{\Bbb Q}$ (in fact we could assume as well that $k$ is a number field or
even only a field of finite type over ${\Bbb Q}$). For any variety $X$ over $%
{\Bbb Q}$, one can choose a model $\cX$ of $X$ over ${\Bbb Z}$, and consider
the number of points $n_{p}(X)$ of the reduction of $\cX$ modulo $p$, for $p$
a prime number. Of course, for some prime numbers $p$, $n_{p}(X)$ may depend
of the model $\cX$, but, for a given $X$, the numbers $n_{p}(X)$ are well
defined for almost all $p$. If we denote by $\cP$ the set of all primes, the
sequence $n_{p}(X)$ is well defined as an element of the ring ${\Bbb Z}^{\cP
}{}^{\prime }$, where, for any ring $R$, we set $R^{\cP}{}^{\prime
}:=\prod_{p\in {\cP}}{R}/\oplus _{p\in {\cP}}{R}$. Moreover, counting points
being additive for disjoint unions and multiplicative for products, the
sequence $n_{p}(X)$ in ${\Bbb Z}^{\cP}{}^{\prime }$ only depends on the
class of $X$ in $\cM$ and may be extended to a ring morphism $n~:\cM
\rightarrow {\Bbb Z}^{\cP}{}^{\prime }$. Setting $n_{p}({\bf L}^{-1})=1/p$,
one may extend uniquely $n$ to a ring morphism $n~:\cM_{{\rm loc}
}\rightarrow {\Bbb Q}^{\cP}{}^{\prime }$. Furthermore, by Lemma 8.1.1 of
\cite{def}, the morphism $n$ factors through the image $\overline{\cM}_{{\rm %
loc}}$ of $\cM_{{\rm loc}}$ in $\widehat{\cM}$.
Recall that we denoted by $%
\overline{\cN}$ the subring of the ring $\widehat{{\cal M}}[[T]]$ which is
generated by the image in $\widehat{{\cal M}}[[T]]$ of ${\cal M}_{{\rm loc}%
}[T]$, $({\bf L}^{i}-1)^{-1}$ and $(1-{\bf L}^{-a}T^{b})^{-1}$ with $i\in
{\Bbb N}$, $a\in {\Bbb N}$, $b\in {\Bbb N}\setminus \{0\}$. Hence, sending $(%
{\bf L}^{i}-1)^{-1}$ and $(1-{\bf L}^{-a}T^{b})^{-1}$ to $({p}^{i}-1)^{-1}$
and $(1-{p}^{-a}T^{b})^{-1}$, respectively, one can uniquely extend $n$ to a
ring morphism $n~:\overline{\cN }\rightarrow {{\Bbb Q}[[T]]}^{\cP}{}^{\prime
}$.

\subsection{}Assume $k={\Bbb Q}$ and that $f_{i}$ and $g_{i}$, for
$i=0,\ldots ,l$, are
polynomials with coefficients in ${\Bbb Q}$. Then, for every prime $p$, we
can consider the $p$-adic cone integral
\[
Z_{{\cal D},p}(s)=\int_{V_{p}}{p}^{-{\rm ord}_{p}(f_{0})s-{\rm ord}
_{p}(g_{0})}d\mu _{p},
\]
with
\[
V_{p}=\left\{ {\bf x}\in {\Bbb Q}_{p}^{m}:{\rm ord}_{p}(f_{i}({\bf x)})\leq
{\rm ord}_{p}(g_{i}({\bf x}))\text{ for }i=1,\ldots ,l\right\} .
\]
Here ${\rm ord}_{p}$ stands for $p$-adic valuation and $d\mu _{p}$ denotes
the additive Haar measure on ${\Bbb Q}_{p}^{m}$, normalized
so that ${\Bbb Z}%
_{p}^{m}$ has volume 1.

By \cite{duSG-Abscissa}, for almost all $p$,
\[
Z_{{\cal D},p}(s)= \sum_{k=0}^{w}(p-1)^{|I_{k}|}{p}^{-m}|E_{I_{k}}^{\circ }(%
{\Bbb F}_p)| \prod_{j\in M_{k}}\frac{{p}^{-(A_{j}s+B_{j})}}{(1-{p}
^{-(A_{j}s+B_{j})})}.
\]
Hence Theorem \ref{explicit expression for motivic cone} has the following
corollary:

\begin{corollary}
\label{np}For almost all primes $p$,
\[
n_{p}(Z_{{\cal D},{\rm geom}}(s))=Z_{{\cal D},p}(s).
\]
(In the left hand side $T={\bf L}^{-s}$, while in the right hand side $T={p}%
^{-s}$.)
\end{corollary}

In particular, we have the following relationship between the motivic zeta
function of a Lie algebra and the zeta functions counting finite index $%
{\Bbb Z}_{p}$-subalgebras:

\begin{corollary}
\label{spelie}Let $L$ be a Lie-algebra of finite rank defined over ${\Bbb Z}$%
. For almost all $p$,
\[
n_{p}(P_{L\otimes {\Bbb Q}[[t]],{\cal X}_{t}^{\leq }}({\bf L}^{-s}))=\zeta
_{L\otimes {\Bbb Z}_{p}}(s).
\]
\end{corollary}

\subsection{}In fact, we should mention that a result such as Corollary
\ref{np} holds in
much wider generality. Let $k$ be a finite extension of $\QQ$ with ring of
integers $\cO$ and $R=\cO
[\frac{1}{N}]$, with $N$ a non zero multiple of the discriminant. For $x$ a
closed point of $\spec R$, we denote by $K_{x}$ the completion of the
localization of $R$ at $x$, by $\cO_{x}$ its ring of integer, and by $\FF_{x}
$ the residue field at $x$, a finite field of cardinality $q_{x}$. Let $f(x)$
be a polynomial in $k[x_{1},\ldots ,x_{m}]$ (or more generally a definable
function in the first order language of valued fields with variables and
values taking place in the valued field and with coefficients in $k$) and
let $\varphi $ be a formula in the language of valued fields with
coefficients in $k$, free variables $x_{1},\ldots ,x_{m}$ running over the
valued field and no other free variables. Now set
$W_{x}:=\{y\in \cO
_{x}^{m} : \varphi (y)\,\text{holds}\}$ and define
\[
I_{\varphi ,f,x}(s)=\int_{W_{x}}|f|_{x}^{s}|dx|_{x},
\]
for $x$ a closed point of ${\rm Spec}\,\cO$. One should remark these
integrals are of much more general nature than the previously considered
cone integrals, since quantifiers are now allowed in $\varphi $. In this
more general setting we still want a geometrical object which specializes,
for almost all $x$ to the local integrals $I_{\varphi ,f,x}(s)$. In \cite
{def} this aim is achieved, at the cost of replacing the ``naive'' ring $\cM%
_{{\rm loc}}$ by a Grothendieck group of Chow motives\footnote{%
For an indication of why Chow motives arise naturally in these questions, see
\ref{isoge}}. More precisely, it is shown in \cite{def}, that there exists a
canonical rational function $I_{\varphi ,f,{\rm mot}}(T)$ with coefficients
in an appropriate Grothendieck ring of Chow motives, such that, for almost
all closed points $x$ in $\spec \cO$, $I_{\varphi ,f,{\rm mot}}(T)$
specializes - after taking the trace of the Frobenius on an {{\'e}}tale
realization and after setting $T=q_{x}^{-s}$ - to the $p$-adic integral $%
I_{\varphi ,f,x}(s)$.

\subsection{}\label{isoge}To what extent are the constructible sets
$E_{I_{k}}^{\circ }$
in Theorem \ref{explicit expression for motivic cone} unique? The left hand
side of (\ref{form}) being canonical, it follows that the class of the right
hand side is independent of the resolution. On the other hand it is unclear
how one can deduce the motivic zeta function from the $p$-adic ones. The
first problem is that distinct varieties over ${\Bbb Q}$ may have the same
number of points in ${\Bbb F}_{p}$, which is in particular the case for
isogenous elliptic curves over ${\Bbb Q}$. Much deeper is the fact, due to
Faltings \cite{Faltings}, that if, for almost all $p$, $n_{p}(E)=n_{p}(E^{%
\prime })$, with $E$ and $E^{\prime }$ elliptic curves over ${\Bbb Q}$, then
$E$ and $E^{\prime }$ are isogenous. Since isogenous elliptic curves define
the same Chow motive, one sees Chow motives appearing in a natural way in
the play. Let us recall (see \cite{Sc} for more details), that for $k$ a
field, a Chow motive over $k$ is just a triple $(S,p,n)$ with $S$ proper and
smooth over $k$, $p$ an idempotent correspondence with coefficients in $%
{\Bbb Q}$ on $X$, and $n$ in ${\Bbb Z}$. One endows Chow motives over $k$
with the structure of a pseudo-abelian category which we denote by ${\rm Mot}%
_{k}$ and we denote by $K_{0}({\rm Mot}_{k})$ its Grothendieck group which
may also be defined as the abelian group associated to the monoid of
isomorphism classes of objects in ${\rm Mot}_{k}$ with respect to the
natural sum $\oplus $. There is also a tensor product on ${\rm Mot}_{k}$
inducing a product on $K_{0}({\rm Mot}_{k})$ which provides $K_{0}({\rm Mot}
_{k})$ with a natural ring structure. Assume now that $k$ is of
characteristic zero. By a result of Gillet and Soul{{\'e}} \cite{G-S} and
Guill{{\'e}}n and Navarro Aznar \cite{G-N} there exists a unique morphism of
rings
\[
\chi _{c}:\cM \longrightarrow K_{0}({\rm Mot}_{k})
\]
such that $\chi _{c}([S])=[h(S)]$ for $S$ projective and smooth, where $h(S)$
denotes the Chow motive associated to $S$, {\it i.e.} $h(S)=(S,{\rm id},0)$.
Let us still denote by $\LL$ the image of $\LL$ by $\chi _{c}$. Since $\LL
=[({\rm Spec}\,k,{\rm id},-1)]$, it is invertible in $K_{0}({\rm Mot}_{k})$,
hence $\chi _{c}$ can be extended uniquely to a ring morphism
\[
\chi _{c}:\cM_{{\rm loc}}\longrightarrow K_{0}({\rm Mot}_{k}).
\]

\begin{problem}\label{pb}Assume $k = {\Bbb Q}$ (one can ask a similar question
when $k$ is a number field, or more generaly of finite type over
${\Bbb Q}$). Is it true that for $\alpha$ in
$\cM_{\rm loc}$ if $n (\alpha) = 0$ then $\chi_c
(\alpha) = 0$?
\end{problem}

A positive answer to this question would imply that the knowledge of the
motivic integrals may be deduced from that of the corresponding $p$-adic
integrals.

\subsection{}Since in the simplest examples, as the ones given in \S \ref
{last}, the zeta functions $\zeta _{L_{p}}(s)$ are just rational functions
of $p$ and $p^{-s}$ and the corresponding motivic zeta functions rational
functions of $\LL$ and $\LL^{-s}$, we give an example showing that non
trivial motives may indeed appear.

Let $L$ be the class two nilpotent Lie algebra over ${\Bbb Z}$ of dimension
9 as a free ${\Bbb Z}$-module given by the following presentation:
\[
L=\left\langle
\begin{array}{c}
x_{1},x_{2},x_{3},x_{4},x_{5},x_{6},y_{1},y_{2},y_{3}:(x_{1},x_{4})=y_{3},(x
_{1},x_{5})=y_{1},(x_{1},x_{6})=y_{2}
\\
(x_{2},x_{4})=y_{2},(x_{2},x_{6})=y_{1},(x_{3},x_{4})=y_{1},(x_{3},x_{5})=y_{3}
\end{array}
\right\rangle
\]
where all other commutators are defined to be 0. This Lie algebra was
discovered by the first author to answer negatively the following question
asked by Grunewald, Segal and Smith in \cite{GSS}: if $L$ is a Lie algebra
over ${\Bbb Z}$, do there exist finitely many rational function $\Phi
_{1}(X,Y),\ldots ,\Phi _{r}(X,Y)$ such that for each prime $p$ there exists $%
i\in \left\{ 1,\ldots ,r\right\} $ and $\zeta _{L_{p}}^{\triangleleft
}(s)=\Phi _{i}(p,p^{-s})?$ Let $E$ be the elliptic curve $Y^{2}=X^{3}-X.$
One can show, by a direct calculation (see \cite{duSG-class2}), that there
exist two non zero rational functions $P_{1}(X,Y)$ and $P_{2}(X,Y)\in {\Bbb Q%
}(X,Y)$ such that, for almost all primes $p$,
\[
\zeta _{L_{p}}^{\triangleleft }(s)=P_{1}(p,p^{-s})+|E({\Bbb F}
_{p})|P_{2}(p,p^{-s}).
\]
By the result of Faltings already alluded to, it follows that the Chow
motive of the curve $E$ is canonically attached to the Lie algebra $L$. It
seems also quite natural to guess, but this has to be checked, that the
calculation of the associated motivic cone integral gives
\[
P_{1}({\bf L},{\bf L}^{-s})+[E]P_{2}({\bf L},{\bf L}^{-s}).
\]

To see where the elliptic curve is hidden in this Lie algebra, consider $
\det ((x_{i},x_{i+3}))=0$ which gives an equation for the projective version
of $E.$

\subsection{}There is one case in which one can deduce an expression for
the motivic zeta
function from the corresponding expression for the $p$-adic zeta functions:

\begin{theorem}
Let $L$ be a Lie algebra over ${\Bbb Z}$ additively isomorphic to ${\Bbb Z}%
^{d}.$ Suppose that there exist two rational functions $\Phi (X,Y)$ and $%
\Psi (X,Y)\in {\Bbb Q}(X,Y)$ such that, for almost all primes $p,$ $\zeta
_{L_{p}}(s)=\Phi (p,p^{-s})$ and $P_{L\otimes {\Bbb Q}[[t]],{\cal X}_{t}^{%
\leq }}({\bf L}^{-s})=\Psi ({\bf L},{\bf L}^{-s}).$ Then $\Phi (X,Y)=\Psi
(X,Y).$
\end{theorem}

\begin{proof} This follows from Corollary \ref{spelie} and the fact that if $%
\Phi (p,p^{-s})=\Psi (p,p^{-s})$ for almost all primes $p$ then $\Phi
(X,Y)=\Psi (X,Y).$ \end{proof}

\section{Topological zeta functions}

\subsection{}In \cite{Denef-Loeser-topological} a topological zeta function was
associated to each polynomial $f\in {\Bbb C}[X]$. The function was defined
from the explicit formula in terms of a resolution of singularities for the
associated local Igusa zeta functions of $f$. It can be interpreted as
considering the explicit formula as $p\mapsto 1.$ The heart of the paper
\cite{Denef-Loeser-topological} was the proof that this definition was
independent of the resolution. The proof given there depended essentially on
using the $p$-adic analysis of these local zeta functions and the associated
explicit formulae. However the current motivic setting gives a proof of the
independence of the resolution without going to the local zeta functions as
we shall explain.

In the setting of groups and Lie algebras it was unclear that such a
topological zeta function could be associated to every nilpotent group or $%
{\Bbb Z}$-Lie algebra. This depended on producing an explicit formula. Since
\cite{duSG-Abscissa} now provides such an explicit formula for zeta
functions of nilpotent groups and Lie algebras and more generally for any
cone integral, we can define the associated topological zeta function and
use the setting of the current paper to prove that it is independent of the
resolution and subdivision of the cone.

We now proceed similarly as in section 2.3 of
\cite{Denef-Loeser-motivic Igusa} to define
topological zeta functions.

Let ${\cal M}[{\bf L}^{-s}]_{{\rm loc}}^{\prime }$ be the subring of ${\cal M%
}[{\bf L}^{-s}]_{{\rm loc}}$ generated by the ring of polynomials ${\cal M}_{%
{\rm loc}}[{\bf L}^{-s}]$ and by the quotients $({\bf L}-1)(1-{\bf L}
^{-Ns-n})^{-1}$ for $(N,n)\in {\Bbb N}^{2}\backslash (0,0)$.

Recall we denoted by
$\overline{\cM}_{{\rm %
loc}}$ the image of $\cM_{{\rm loc}}$ in $\widehat{\cM}$. In
 particular $\overline{\cM}_{{\rm %
loc}}$
is
a quotient without $({\bf L}-1)$-torsion of
${\cal M}_{{\rm loc}}$.

By expanding ${\bf L}^{-s}$ and $({\bf L}-1)(1-{\bf L}^{-Ns-n})^{-1}$ into
series in ${\bf L}-1,$ one gets a canonical morphism of algebras
\[
\varphi :{\cal M}[{\bf L}^{-s}]_{{\rm loc}}^{\prime }\rightarrow \overline{%
{\cal M}}_{{\rm loc}}[s][(Ns+n)^{-1}]_{(N,n)\in {\Bbb N}^{2}\backslash
(0,0)}[[{\bf L}-1]],
\]
where $[[{\bf L}-1]]$ denotes completion with respect to the ideal generated
by ${\bf L}-1$.
Taking the quotient of $\overline{{\cal M}}_{{\rm loc}}[s][(Ns+n)^{-1}
]_{(N,n)\in {\Bbb N}^{2}\backslash (0,0)}[[{\bf L}-1]]$ by the ideal
generated by ${\bf L}-1,$ one obtains the evaluation morphism:
\begin{multline*}
{\rm ev}_{{\bf L}=1}:\overline{{\cal M}}_{{\rm loc}}[s][(Ns+n)^{-1}%
]_{(N,n)\in {\Bbb N}^{2}\backslash (0,0)}[[{\bf L}-1]]\longrightarrow \\
\left( \overline{{\cal M}}_{{\rm loc}}/{\bf L}-1\right)
[s][(Ns+n)^{-1}]_{(N,n)\in {\Bbb N}^{2}\backslash (0,0)}.
\end{multline*}

We denote by $\chi_{{\rm top}}(X)$ the usual Euler characteristic of a
smooth, projective $k$-scheme (say in {\'e}tale $\overline{{\Bbb Q}_{l}}$
-cohomology). By \S\kern .15em
6.1 of
\cite{Denef-Loeser-Arcs}, this factorises through a morphism $\chi_{{\rm top}}:
\overline{{\cal M}}_{{\rm loc}}\rightarrow {\Bbb N}$. Since $\chi_{{\rm
top}}({\bf L})=1$, this induces then a
morphism
\[
\chi_{{\rm top}}:\left( \overline{{\cal M}}_{{\rm loc}}/{\bf L}-1\right)
[s][(Ns+n)^{-1}]_{(N,n)\in {\Bbb N}^{2}\backslash (0,0)}\rightarrow {\Bbb Z}%
[s][(Ns+n)^{-1}]_{(N,n)\in {\Bbb N}^{2}\backslash (0,0)}.
\]

\begin{definition}
The {\em topological cone zeta function} $Z_{{\cal D},{\rm top}}(s)$
associated to the cone data ${\cal D}$ is defined to be
\[
Z_{{\cal D},{\rm top}}(s)
:=\left( \chi _{{\rm top}}\circ {\rm ev}_{{\bf L}%
=1}\circ \varphi \right) \left( Z_{{\cal D},{\rm geom}}(s)\right) .
\]
\end{definition}

\begin{definition}Let
$L$ be a $k[[t]]$-algebra additively isomorphic to $\left( k[[t]]\right)
^{d}$ or a ${\Bbb Z}$-Lie algebra additively isomorphic to ${\Bbb Z}^{d}.$
Let $L$ be an algebra additively isomorphic to ${\Bbb Z}^{d}.$
We define the
{\em topological zeta function }$Z_{L,{\rm top}}(s)$ of $L$ to be
\[
Z_{{L},{\rm top}}(s)
:=\left( \chi _{{\rm top}}\circ {\rm ev}_{{\bf L}%
=1}\circ \varphi \right) \left( Z_{{L},{\rm geom}}(s)\right) .
\]
\end{definition}

These definitions make sense since $Z_{{\cal D},{\rm geom}}(s)$
and
$Z_{{L},{\rm geom}}(s)$
belong to
${\cal M}[{\bf L}^{-s}]_{{\rm loc}}^{\prime }$
by Theorem \ref{explicit expression for motivic cone}.

We can get an explicit formula for the topological zeta functions from that
for the associated motivic zeta function. We use the notation defined prior
to the statement of Theorem \ref{explicit expression for motivic cone}. We
define the subset $W_{{\rm top}}\subset \left\{ 0,\ldots,w\right\}$ such
that $k\in W_{{\rm top}}$ if and only if $|I_{k}|=|M_{k}|.$ These are the
indices of pieces of the cone whose dimension is equal to the dimension of
the vector space spanned by all those basis elements with non-zero
coefficients somewhere in the piece.

\begin{proposition}
Let ${\cal D}$ be a set of cone integral data. Then with the notation above:
\[
Z_{{\cal D},{\rm top}}(s)=\sum_{k\in W_{{\rm top}}}\chi _{{\rm top}%
}(E_{I_{k}}^{\circ })\prod_{j\in M_{k}}\frac{1}{(A_{j}s+B_{j})}.
\]
\end{proposition}

Some pieces $R_{j}$ of the cone may be missing in this expression.
Note that in the special case of the Igusa zeta function we do see all the
pieces $R_{j}$ since the cone there is the full positive quadrant.
Nevertheless, in the more general setting of a cone integral, the constants $%
A_{j}$ and $B_{j}$ ($j\in \left\{ 1,\ldots ,q\right\} $) will appear in this
formula somewhere whenever the one-dimensional edge $R_{j}$  sits on an open
simplicial piece of the cone of dimension corresponding to the number of
non-zero entries of $R_{j}$.

\section{Examples}

\label{last}

\subsection{Free abelian}

Let $L={\Bbb Z}^{d}.$ Then
\begin{eqnarray*}
P_{L\otimes k[[t]], \cX_t^{\leq}}(s) &=&(1-{\bf L}^{-s})^{-1}\ldots(1-{\bf L}
^{-s+d-1})^{-1} \\
Z_{L,{\rm geom}}(s) &=&\frac{(1-{\bf L}^{-1})}{(1-{\bf L}^{-s-d})}\ldots%
\frac{(1-{\bf L}^{-d})}{(1-{\bf L}^{-s-1})} \\
Z_{L,{\rm top}}(s) &=&\frac{1}{(s+d)\ldots(s+1)}.
\end{eqnarray*}

The calculation of the integral in the $p$-adic setting translates
immediately to a calculation of the motivic integral, proving the above
identities.

\subsection{Heisenberg}

Let
\[
L=\left(
\begin{array}{lll}
0 & {\Bbb Z} & {\Bbb Z} \\
& 0 & {\Bbb Z} \\
&  & 0
\end{array}
\right) .
\]
Then
\begin{eqnarray*}
P_{L\otimes k[[t]], \cX_t^{\leq}}(s) &=&(1-{\bf L}^{-s})^{-1}(1-{\bf
L}^{1-s})^{-1}(1-{\bf %
L}^{2-2s})^{-1}(1-{\bf L}^{3-2s})^{-1}(1-{\bf L} ^{3-3s}) \\
Z_{L,{\rm geom}}(s) &=&\frac{(1-{\bf L}^{-1})(1-{\bf L}^{-2})(1-{\bf L}%
^{-3})(1-{\bf L}^{-3s-6})}{(1-{\bf L}^{-s-3})(1-{\bf L} ^{-s-2})(1-{\bf L}
^{-2s-4})(1-{\bf L}^{-2s-3})} \\
Z_{L,{\rm top}}(s) &=&\frac{3}{2\left( s+3\right) \left( s+2\right) \left(
2s+3\right) }.
\end{eqnarray*}

Using Proposition \ref{conversion to triangular},
\[
Z_{L,{\rm geom}}(s)=\frac{(1-{\bf L}^{-1})(1-{\bf L}^{-2})(1-{\bf L}^{-3})}{%
(1-{\bf L}^{-1})^{3}}\int_{A}{\bf L}^{-{\rm ord}_{t}(x)(s+2)}{\bf L}^{-{\rm %
ord}_{t}(y)(s+1)}{\bf L}^{-{\rm ord}_{t}(z)s}d\mu
\]
where
\[
A=\left\{ \left( x,y,z\right) \in k[[t]]^{3}:{\rm ord}_{t}(xy)\geq {\rm ord}
_{t}(z)\right\} .
\]

The polynomial associated to the cone integral data already defines a
non-singular space consisting only of normal crossings. So no resolution of
singularities is required. The associated cone $\overline{D_{T}}$ is the
span of the four one-dimensional pieces where we record the associated $%
A_{j} $ and $B_{j}$:
\begin{eqnarray*}
R_{1} &=&(1,0,1),A_{1}=2,B_{1}=4 \\
R_{2} &=&(1,0,0),A_{2}=1,B_{2}=3 \\
R_{3} &=&(1,0,1),A_{3}=2,B_{3}=3 \\
R_{4} &=&(0,1,0),A_{4}=1,B_{4}=2.
\end{eqnarray*}

The following is a table of the pieces $R_{k}$ $k\in \{0,\ldots,w\}$
corresponding to a suitable decomposition of this cone where $%
R_{i_{1}}\ldots R_{i_{r}}$ denotes the open piece spanned by $%
R_{i_{1}},\ldots,R_{i_{r}}$:

\[
\begin{array}{llll}
R_{k} & |I_{k}| & |M_{k}| & E_{I_{k}}^{\circ } \\
0 & 0 & 0 & \left( {\bf L}-1\right) ^{3} \\
R_{1} & 2 & 1 & \left( {\bf L}-1\right)  \\
R_{2} & 1 & 1 & \left( {\bf L}-1\right) ^{2} \\
R_{3} & 2 & 1 & \left( {\bf L}-1\right)  \\
R_{4} & 1 & 1 & \left( {\bf L}-1\right) ^{2} \\
R_{1}R_{2} & 2 & 2 & \left( {\bf L}-1\right)  \\
R_{1}R_{3} & 3 & 2 & 1 \\
R_{2}R_{3} & 3 & 2 & 1 \\
R_{2}R_{4} & 2 & 2 & \left( {\bf L}-1\right)  \\
R_{3}R_{4} & 2 & 2 & \left( {\bf L}-1\right)  \\
R_{1}R_{2}R_{3} & 3 & 3 & 1 \\
R_{2}R_{3}R_{4} & 3 & 3 & 1
\end{array}
\]
Using our explicit formula it is possible then to derive the expression for $%
Z_{L,{\rm geom}}(s)$ given above.

To derive the expression for $Z_{L,{\rm top}}(s)$ note that $\chi_{{\rm top}
}(\left( {\bf L}-1\right) ^{i})=0$ for $i>0.$ Hence those pieces $R_{k}$ for
which $|I_{k}|=|M_{k}|$ and $\chi_{{\rm top}}(E_{I_{k}}^{\circ })\neq 0$ are
$R_{1}R_{2}R_{3}$ and $R_{2}R_{3}R_{4}.$ This gives the following expression
for $Z_{L,{\rm top}}(s):$%
\begin{eqnarray*}
Z_{L,{\rm top}}(s) &=&\frac{1}{(2s+4)(s+3)(2s+3)}+\frac{1}{(s+3)(2s+3)(s+2)}
\\
&=&\frac{3}{2\left( s+3\right) \left( 2s+3\right) \left( s+2\right) }.
\end{eqnarray*}
Of course one could directly deduce this by applying the function $\chi_{%
{\rm top}}\circ {\rm ev}_{{\bf L}=1}\circ \varphi $ to the final expression
for $Z_{L,{\rm geom}}(s)$. However it is instructive to see the explicit
formula for $Z_{L,{\rm top}}(s)$ at work.

\subsection{${\frak sl}_{2}$}

Let
\[
L=\left\{ \left(
\begin{array}{ll}
a & b \\
c & -a
\end{array}
\right) :(a,b,c)\in {\Bbb Z}^{3}\right\} .
\]
Then
\begin{eqnarray*}
P_{L\otimes k[[t]], \cX_t^{\leq}}(s) &=&(1-{\bf L}^{-s})^{-1}(1-{\bf
L}^{1-s})^{-1}(1-{\bf %
L}^{2-2s})^{-1}(1-{\bf L}^{1-2s})^{-1}(1-{\bf L} ^{1-3s}) \\
Z_{L,{\rm geom}}(s) &=&\frac{(1-{\bf L}^{-1})(1-{\bf L}^{-2})(1-{\bf L}
^{-3})(1-{\bf L}^{-3s-8})}{(1-{\bf L}^{-s-3})(1-{\bf L}^{-s-2})(1-{\bf L}%
^{-2s-4})(1-{\bf L}^{-2s-5})} \\
Z_{L,{\rm top}}(s) &=&\frac{3s+8}{2\left( s+3\right) \left( s+2\right)
^{2}\left( 2s+5\right) }.
\end{eqnarray*}

Using Proposition \ref{conversion to triangular} we have the following
description of $Z_{L,{\rm geom}}(s)$%
\[
Z_{L,{\rm geom}}(s)=\frac{(1-{\bf L}^{-1})(1-{\bf L}^{-2})(1-{\bf L}^{-3})}{
(1-{\bf L}^{-1})^{3}}\int_{W}{\bf L}^{-{\rm ord}_{t}(a)(s+2)}{\bf L}^{-{\rm %
ord}_{t}(x)(s+1)}{\bf L}^{-{\rm ord}_{t}(z)s}d\mu
\]
where
\[
W=\left\{ \left(
\begin{array}{lll}
a & b & c \\
0 & x & y \\
0 & 0 & z
\end{array}
\right) \in {\rm Tr}_{3}(R):\left.
\begin{array}{c}
v(x)\leq v(4zb) \\
v(x)\leq v(4by) \\
v(zx)\leq v(ax^{2}+4cxy-4by^{2})
\end{array}
\right. \right\} .
\]

The polynomial associated to the cone integral data is singular in this
case. There have been several calculations of the zeta function of ${\frak sl%
}_{2}({\Bbb Z}_{p})$ for $p>2$ (see \cite{Ilani-sl2} and \cite{dus-sl2})
which give the answer above with ${\bf L}$ replaced by $p$. However, it is
not  possible to deduce directly that the above expression is a correct
formula for the corresponding motivic zeta function, cf. the discussion at
the end of Section \ref{Section:Motivic cone integrals}. However, recently,
the first author and Ph.D. student Gareth Taylor \cite{duS and Taylor}, have
used the approach of applying resolution of singularities to the
non-singular polynomial associated to the cone integral data of the integral
above to deduce the expression for the $p$-adic integrals. The advantage of
this approach is that it formally translates into a calculation of the
corresponding motivic zeta function. Alternatively, one can observe that the
resolution of the polynomial gives rise to $E_{I}$ which are generated out
of ${\bf L}$, {\it i.e.} we don't get any exotic varieties, so that there
exists a rational function $\Psi (X,Y)\in {\Bbb Q}(X,Y)$ such that the
motivic zeta function is $\Psi ({\bf L},{\bf L}^{-s})$. Hence, by Corollary
\ref{spelie}, one is able to deduce that the rational function describing
the $p$-adic expressions is the same as $\Psi (X,Y).$


\begin{thebibliography}{99}
\bibitem{BL}  A. Beauville and Y. Laszlo, {\it Conformal blocks and
generalized theta functions}, Comm. Math. Phys. {\bf 164} (1994), 385--419.

\bibitem{Bourbaki-Lie algebra}  N. Bourbaki, {\it Lie Groups and Lie Algebras%
}, Part I, (1989) Springer-Verlag, Berlin.

\bibitem{Danilov}  V.I. Danilov, {\it The geometry of toric varieties},
Russian Math. Surveys {\bf 33} (1978), 97--154.

\bibitem{rat}  J. Denef, {\it The rationality of the Poincar{{\'e}} series
associated to the $p$-adic points on a variety}, Invent. Math., {\bf 77}
(1984), 1--23.

\bibitem{Denef-Loeser-topological}  J. Denef and F. Loeser, {\it Caract{{\'e}%
}ristiques d'Euler-Poincar{{\'e}}, fonctions z{{\^e}}ta locales et
modifications analytiques}, J. of the Amer. Math. Soc. {\bf 5} (1992),
705--720.

\bibitem{Denef-Loeser-motivic Igusa}  J. Denef and F. Loeser, {\it Motivic
Igusa zeta functions}, J. Algebraic Geom., {\bf 7} (1998), 505--537.

\bibitem{Denef-Loeser-Arcs}  J. Denef and F. Loeser, {\it Germs of arcs on
singular algebraic varieties and motivic integration}, Invent. Math., {\bf %
135} (1999), 201--232.

\bibitem{MK}  J. Denef and F. Loeser, {\it Motivic integration, quotient
singularities and the McKay correspondence}, preprint February 1999.

\bibitem{def}  J. Denef and F. Loeser, {\it Definable sets, motives and $p$%
-adic integrals}, preprint October 1999.

\bibitem{duS-Annals}  M.P.F. du Sautoy, {\it Finitely generated groups,
$p$-adic analytic groups and Poincar{{\'e}} series}, Annals of Math. {\bf
137}
(1993), 639--670.

\bibitem{dus-sl2}  M.P.F. du Sautoy, {\it The zeta function of ${\frak sl}%
_{2}({\Bbb Z}),$} to appear in Forum Mathematicum.

\bibitem{duS and Taylor}  M.P.F. du Sautoy and G. Taylor, {\it The zeta
function of ${\frak sl}_{2}$ and resolution of singularities}, in
preparation.

\bibitem{duSG-Abscissa}  M.P.F. du Sautoy and F.J. Grunewald, {\it Analytic
properties of zeta functions and subgroup growth}, M.P.I. preprint series
1999 (87).

\bibitem{duSG-class2}  M.P.F. du Sautoy and F.J. Grunewald, {\it Counting
subgroups of finite index in nilpotent groups of class 2}, in preparation.

\bibitem{Faltings}  G. Faltings, {\it Endlichkeitss{{\"a}}tze f{{\"u}}r
abelsche {V}ariet{{\"a}}ten {{\"u}}ber {Z}ahlk{{\"o}}rpern}, Invent. Math.
{\bf 73} (1983), 349--366.

\bibitem{G-S}  H. Gillet, C. Soul{{\'e}}, {\it Descent, motives and
$K$-theory}, J. Reine Angew. Math., {\bf 478} (1996), 127--176.

\bibitem{GSS}  F.J. Grunewald, D. Segal and G.C. Smith, {\it Subgroups of
finite index in nilpotent groups}, Invent. Math., {\bf 93} (1988), 185--223.

\bibitem{G-N}  F. Guill{{\'e}}n, V. Navarro Aznar, {\it Un crit{{\`e}}re
d'extension d'un foncteur d{{\'e}}fini sur les sch{{\'e}}mas lisses},
preprint (1995), revised (1996).

\bibitem{Hall}  P. Hall, {\it Nilpotent Groups}, Queen Mary College Math.
Notes, (1969).

\bibitem{Ilani-sl2}  I. Ilani, {\it Zeta functions related to the group $%
{\rm SL}_{2}({\Bbb Z}_{p})$}, Israel J. of Math., {\bf 109} (1999), 157--172.

\bibitem{Kac}  V. Kac, {\it Infinite dimensional Lie algebras}, third
edition, Cambridge University Press (1990).

\bibitem{Mathieu}  O. Mathieu, {\it Classification of simple graded Lie
algebras of finite growth}, Invent. Math., {\bf 108} (1992) 455--519.

\bibitem{Pressley-Segal}  A. Pressley and G. Segal, {\it Loop Groups},
Oxford Mathematical Monographs, Oxford University Press (1986).

\bibitem{Reutenauer}  C. Reutenauer, {\it Free Lie Algebras}, London Math.
Soc. Monographs New Series 7, Clarendon Press, Oxford (1993).

\bibitem{Sc}  A. Scholl, {\it Classical motives}, in {\sl Motives}, U.
Jannsen, S. Kleiman, J.-P. Serre Ed., Proceedings of Symposia in Pure
Mathematics, Volume 55 Part 1 (1994), 163--187.

\bibitem{Segal-polycyclic groups}  D. Segal, {\it Polycyclic Groups},
Cambridge Tracts in Mathematics {\bf 82}, Cambridge University Press (1983).

\bibitem{Segal-Wilson}  G. Segal and G. Wilson, {\it Loop groups and
equations of KdV type}, Pub. Math. I.H.E.S., {\bf 61 }(1985), 5--65.

\bibitem{Shalev-New Horizons}  A. Shalev, {\it Lie methods in the theory of
pro-$p$ groups}, in New Horizons in Pro-$p$ Groups, edited by
M.P.F. du Sautoy, D. Segal and A. Shalev, Progress in Mathematics {\bf 184}
Birkh{\"a}user (2000).


\bibitem{Strade and Farnsteiner}  H. Strade and R. Farnsteiner, {\it
Modular Lie
algebras and their representations}, Monographs and Textbooks in Pure
and Applied Mathematics {\bf 11}
Marcel Dekker, Inc., New York (1988).



\end{thebibliography}
\end{document}